\documentclass[11pt,english]{article}
\usepackage[sc]{mathpazo}
\usepackage{geometry}
\usepackage{color}
\usepackage{array}
\usepackage{bbding}
\usepackage{multirow}
\usepackage[fleqn]{amsmath}
\usepackage{amssymb}
\usepackage{amsthm}
\usepackage{graphicx}
\usepackage{natbib}
\usepackage{lscape}
\usepackage[hyphens]{url}
\usepackage[hidelinks]{hyperref}
\usepackage{comment}
\usepackage{bm}
\usepackage[title]{appendix}
\usepackage{longtable}
\usepackage{mathtools}
\usepackage{chngcntr}
\usepackage{authblk}
\usepackage{babel}
\usepackage{dsfont}
\usepackage{subcaption}
\usepackage{textgreek}
\usepackage{cleveref}

\makeatletter
\allowdisplaybreaks
\usepackage[mathlines, pagewise]{lineno}
\usepackage{etoolbox} %% <- for \pretocmd, \apptocmd and \patchcmd

\newcommand*\linenomathpatchAMS[1]{%
	\expandafter\pretocmd\csname #1\endcsname {\linenomathAMS}{}{}%
	\expandafter\pretocmd\csname #1*\endcsname{\linenomathAMS}{}{}%
	\expandafter\apptocmd\csname end#1\endcsname {\endlinenomath}{}{}%
	\expandafter\apptocmd\csname end#1*\endcsname{\endlinenomath}{}{}%
}

\expandafter\ifx\linenomath\linenomathWithnumbers
\let\linenomathAMS\linenomathWithnumbers
\patchcmd\linenomathAMS{\advance\postdisplaypenalty\linenopenalty}{}{}{}
\else
\let\linenomathAMS\linenomathNonumbers
\fi

\linenomathpatchAMS{gather}
\linenomathpatchAMS{multline}
\linenomathpatchAMS{align}
\linenomathpatchAMS{alignat}
\linenomathpatchAMS{flalign}

\newtheorem{prop}{Proposition}

\newcommand{\E}{\mathop{{}\mathbb{E}}}
\newcommand{\Var}{\mathop{{}\mathrm{Var}}}

\def\d{\mathop{\rm \!{d}}}

\newcommand{\probP}{\text{I\kern-0.15em P}}

\usepackage{array}

\usepackage{algorithmicx}
\usepackage{algpseudocode}
\usepackage{xcolor}
\usepackage{mdframed}

\definecolor{myboxcolor}{RGB}{220,220,220}

\newmdenv[
  backgroundcolor=myboxcolor,
  roundcorner=10pt,
  innertopmargin=0pt,
  innerbottommargin=6pt,
  linewidth=1pt, % Add a 1pt black border
  linecolor=black % Set the border color to black
]{coloredbox}

\begin{document}
	
\title{Analysis of Transshipment in Three-Sided Meal Delivery Services via Microhubs}
\author[1]{Linxuan Shi}
\author[1]{Zhengtian Xu\footnote{Corresponding author. E-mail address: \textcolor{blue}{zhengtian@gwu.edu} (Z. Xu).}}
\affil[1]{\small\emph{Department of Civil and Environmental Engineering, The George Washington University, Washington DC, United States}\normalsize}

%\author[]{}
\date{\today}
\maketitle

\begin{abstract}
\noindent This paper introduces and analyzes a novel transshipment strategy for meal delivery. In this approach, the service area is partitioned into smaller sub-areas, with deliverers assigned to operate exclusively within these sub-areas. Meanwhile, a centrally located microhub functions as a logistic depot to facilitate the batching and transfer of meal packages toward different sub-areas. We model the meal delivery system with transshipment using networked G/G/m queues and analytically approximate two critical system performance metrics---customer waiting time and vehicle miles traveled---to evaluate the effectiveness of the proposed strategy. The performance achieved by transshipment is benchmarked against that of the classic pickup-and-delivery strategy without transshipment, both predicted using continuous approximations. For the latter, we enhance the existing modeling by incorporating the delivery distance profiles of individual orders to better match the meal delivery context. Our comparisons indicate that meal delivery via transshipment outperforms the non-transshipping counterpart across both metrics under either high-demand or low-supply conditions, with particular advantages in servicing larger areas or handling long-distance orders. This conclusion is corroborated by a numerical experiment using empirical meal delivery data from Meituan, which suggests that an optimally configured transshipment strategy can significantly improve service performance for both customers and deliverers during peak lunch hours and in the busiest districts. While transshipment continues to reduce vehicle miles traveled by deliverers during non-peak hours, it results in longer customer waiting times compared to the benchmark without transshipment as demand decreases.
\end{abstract}


\newpage
\section{Introduction}

Meal delivery services connect customers and restaurants through mobile applications, allowing users to order meals online with just a few taps on their smartphones. The platform assigns meal orders to a team of affiliated deliverers, who complete the fulfillment process by picking up the meal package from the restaurant and delivering it to the customers. The lasting COVID-19 pandemic prompted an irreversible shift among consumers from indoor dining to online meal delivery \citep{kaczmarski2024which}, under which takeout channels have become increasingly vital to restaurant operations in the post-pandemic era \citep{shi2024dine}. The significant surge in customer demand has heightened the need for meal delivery platforms to innovate supply management to ensure timely delivery. To prioritize customer satisfaction, platforms have been continuously seeking faster last-mile delivery methods, such as strategic order batching and assignment \citep{liu2021time, simoni2023crowdsourced, ye2024modeling}, advanced vehicle routing strategies \citep{steever2019dynamic, kohar2021capacitated}, and new delivery modes utilizing autonomous vehicles and drones \citep{liu2019optimization, shi2022integer}. As intermediaries in a three-sided service market involving customers, deliverers, and restaurants, meal delivery platforms must carefully consider and manage the differing interests of these stakeholders when introducing new delivery operations \citep{agatz2024transportation}.

\begin{figure}[!t]
	\begin{subfigure}[b]{.5\textwidth}
		\centering
		\includegraphics[width=0.9\linewidth]{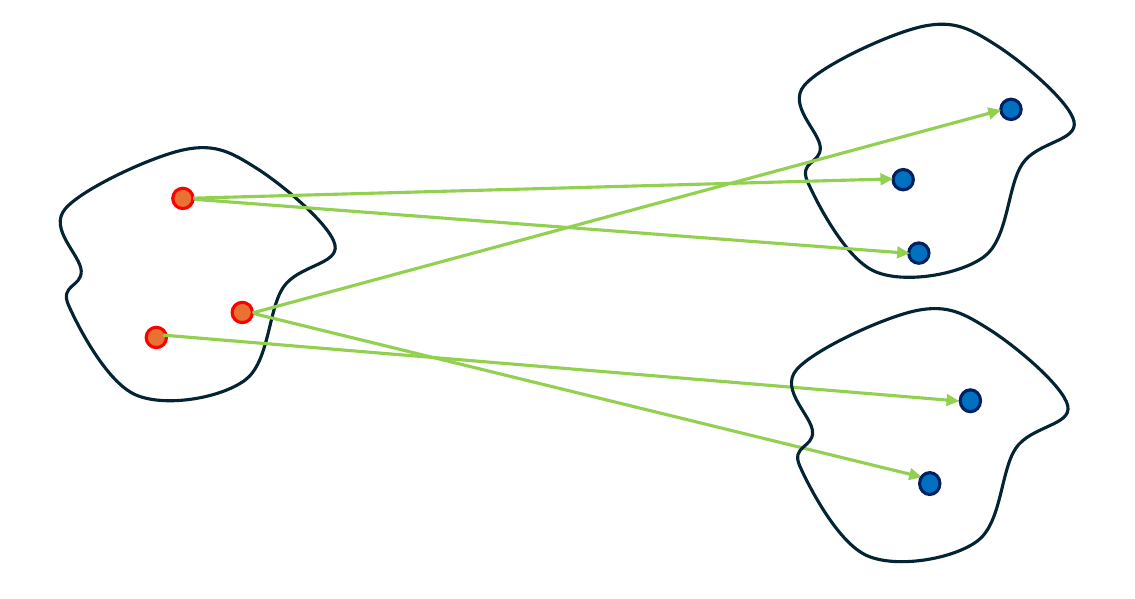}
		\caption{Line-haul logistics without transshipment}
        \label{fig_intro_a}
	\end{subfigure}\hfill
	\begin{subfigure}[b]{.5\textwidth}
		\centering
		\includegraphics[width=0.9\linewidth]{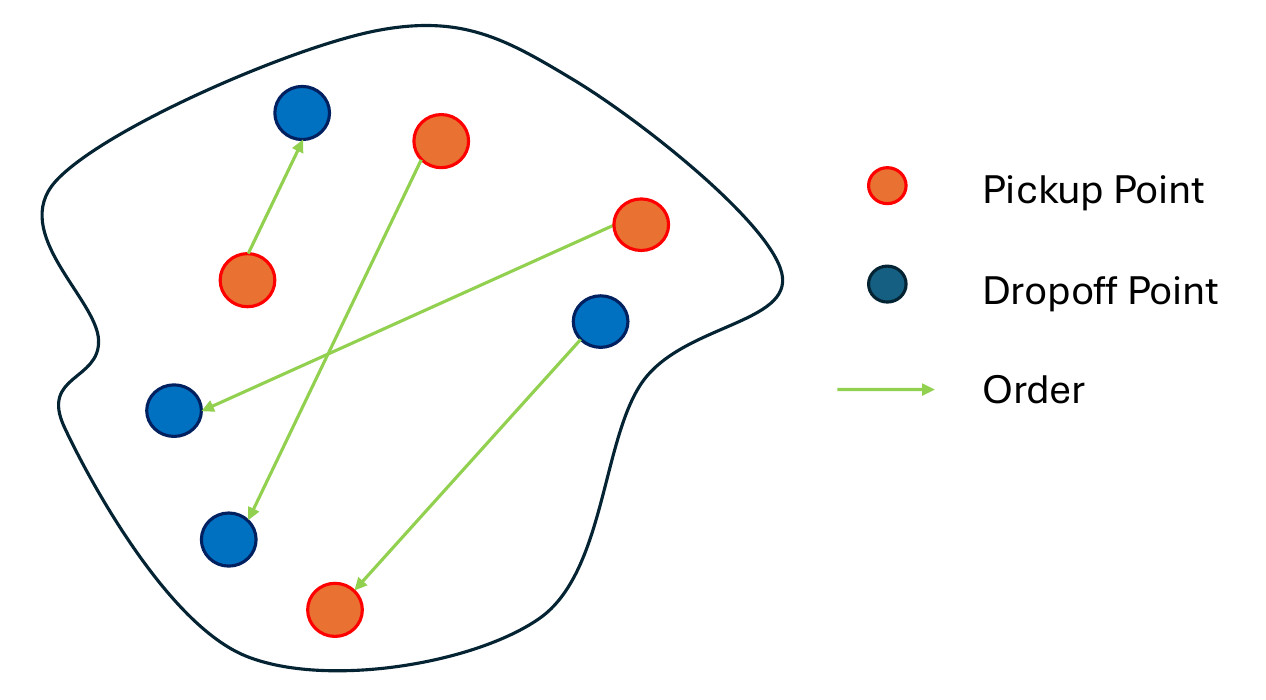}
		\caption{Meal delivery without transshipment}
        \label{fig_intro_b}
	\end{subfigure}
	\begin{subfigure}[b]{.5\textwidth}
		\centering
		\includegraphics[width=0.85\linewidth]{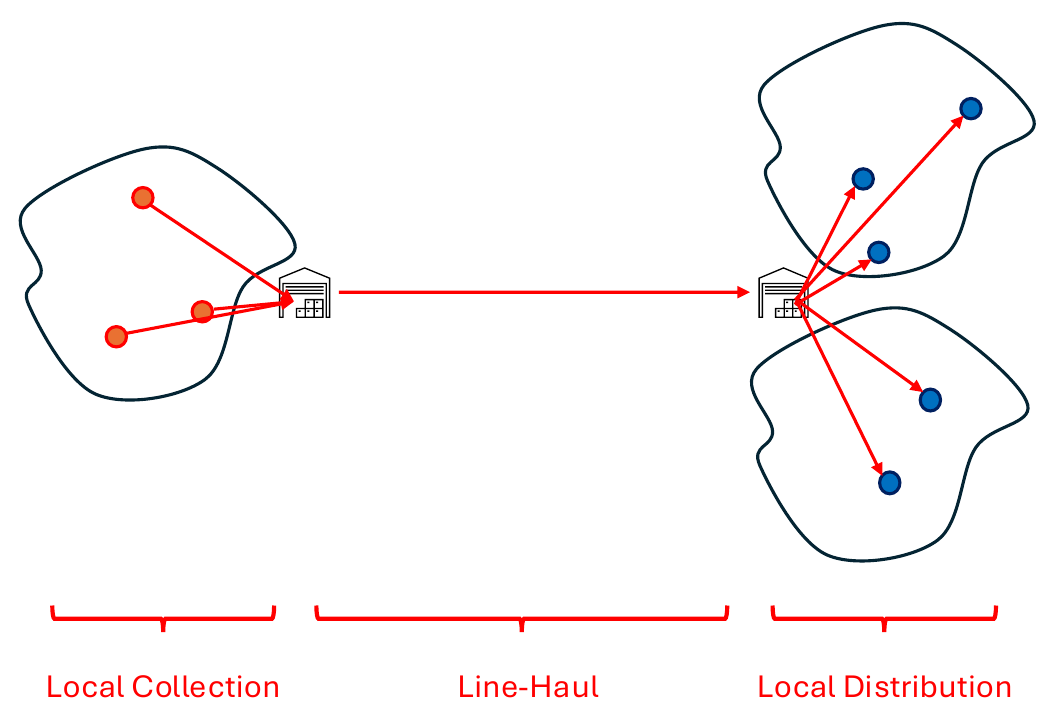}
		\caption{Line-haul logistics with transshipment}
        \label{fig_intro_c}
	\end{subfigure}\hfill
	\begin{subfigure}[b]{.5\textwidth}
		\centering
		\includegraphics[width=0.9\linewidth, height=0.55\linewidth]{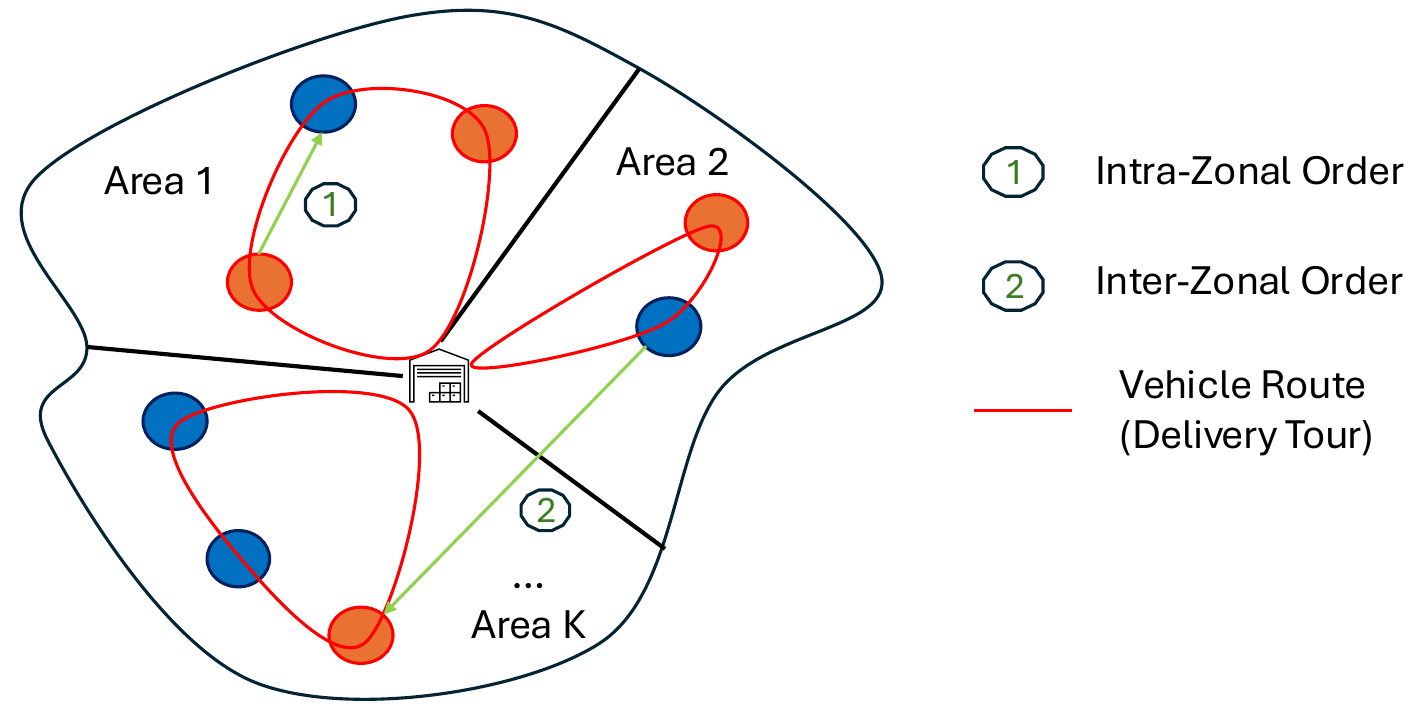}
		\caption{Meal delivery with transshipment}
        \label{fig_intro_d}
	\end{subfigure}
	\caption{Comparison between transshipment in line-haul logistics and meal delivery.}
	\label{fig_intro}
\end{figure}

This study conceptualizes a new operational model for meal delivery by enabling transshipment through a microhub, aiming to enhance overall efficiency and provide long-term benefits to all stakeholders. The concept of microhubs mirrors logistic depots in line-haul logistics, allowing for sorting, storing, and transferring packages, but with a focus on urban neighborhoods. As illustrated in Figure \ref{fig_intro_c}, the availability of these depots re-structures freight distribution into a transshipment-based operation, under which packages are transferred from line-haul trucks to local distribution networks for final delivery. Compared to the delivery model without transshipment, shown in Figure \ref{fig_intro_a}, this approach brings services closer to end customers and enables more flexible operations within targeted service areas. In some use cases, it alleviates congestion by consolidating freight vehicle trips and lowers emissions by facilitating transfers to low- or zero-emission fleets for last-mile deliveries \citep{gunes2024seattle}. 

Unlike the traditional line-haul distribution system shown in Figure \ref{fig_intro_a}, which typically involves one cluster of origins and another of destinations, the meal delivery system features orders with pickup and drop-off points dispersed across the entire region, as illustrated in Figure \ref{fig_intro_b}. In our proposed design, illustrated in Figure \ref{fig_intro_d}, the service region is divided into pie-shaped sub-areas, with one microhub located at the junction, to facilitate localized batching for order collection and distribution. The microhub serves as an intermediate node for the temporal storage and transfer of meal packages, while also coordinating pickup and drop-off tasks both within and across sub-areas. Individual meal packages are picked up from the restaurants where they originated and transported to the microhub. Subsequently, one deliverer collects the package from the microhub and delivers it to the customer. Delivery vehicles are dispatched within each area to perform pickup, drop-off, and transfer tasks, with each vehicle dedicated to serving a specific sub-area. In practical implementation, deliverers possess better knowledge in routing for regions they are familiar with and prefer to operate exclusively within these areas for greater efficiency. In contrast to transshipment in line-haul logistics, which exploits transportation economies of scale by moving goods in large batches over shared long-haul segments, our design of transshipment through a microhub, tailored for intra-city many-to-many distribution contexts, aims to gain efficiency through localized order batching. To evaluate the comparative efficiency gains provided by introducing transshipment, this study employs continuous approximations (CA) to model meal delivery operations. We estimate two key metrics---deliverers' Vehicle Miles Traveled (VMT) and customers' Average Waiting Time---both with and without transshipment.

Compared to traditional optimization-based models used to address the vehicle routing problem with pickup and delivery (VRPPD) at the operational level \citep{toth2002vehicle}, CA models are extensively applied to tackle big-picture questions about logistics and transportation systems \citep{ansari2018advancements}. These applications include evaluating new designs in public transit \citep{chen2017analysis, chen2017connecting} and comparing new strategies with existing ones for dial-a-ride services \citep{xu2020supply, ouyang2023measurement}. While optimization-based approaches have significantly improved in solution quality over the past decade, their heavy formulations may fall short in providing a comprehensive understanding of the problem's inherent characteristics, which could be critical for managerial insights. Additionally, concerns about data quality, particularly regarding precision and uncertainty, cast doubt on the reliability of optimal solutions in practical applications \citep{daganzo1987increasing}. CA models, serving as a complement to optimization models, have proven effective in addressing these challenges. According to \cite{daganzo2012potential}, effective CA models can yield five major benefits: they require fewer data inputs, exhibit reduced computational complexity, provide enhanced system representation, improve transparency, and facilitate deeper insights.

For adopting the CA approach, the challenge lies in describing the system of interest as simply as possible while maintaining accurate estimations, a principle known as creating a parsimonious model. For a meal delivery system, its market conditions can be characterized primarily by three exogenous components: the coverage of service areas, the level of customer demand, and the service capacity provided by deliverers. Given specific market conditions, the platform will make strategic decisions on how to divide the area into smaller sub-areas as well as the strategy for batching and fulfilling pickup and drop-off requests within each sub-area. The system performances are influenced by both the exogenous market conditions and the operating strategy taken by the platform in fulfilling meal deliveries. With the help of CA models, we can address strategic-level questions to identify the market conditions under which introducing transshipment can be advantageous for meal delivery compared to the standard VRPPD strategy.

The remainder of this paper is organized as follows. \Cref{sec:liter_review} first reviews related work to position the contributions of our study to the literature. \Cref{sec:microhub} introduces the design of our proposed strategy for meal delivery with transshipment, followed by a detailed explanation of the derivations to estimate two key system performance metrics: vehicle miles traveled and customer waiting time. Then, \Cref{sec:VRPPD} derives an approximation for VRPPD, tailored to the specific characteristics of meal delivery, which serves as a benchmark for comparison. Comprehensive numerical experiments are conducted in \Cref{sec:numerical} to validate the accuracy of our approximations and to compare the performance of the proposed transshipment strategy against the prevalent VRPPD approach. Finally, \Cref{sec:conclusion} summarizes the findings and concludes the paper.

\section{Literature Review}\label{sec:liter_review}

In the design of freight distribution systems, there are two main classes based on the presence of transshipment: systems with transshipment and those without. The proposed integrated meal delivery system with a microhub falls into the former category. On the other hand, within each partitioned sub-area, the operations of deliverers align with pickup-and-delivery problems (PDPs) without transshipment.

\subsection{Distribution system without transshipment}

Most studies on distribution systems without transshipment focus on identifying the most cost-effective routes to service demand nodes within a network, such as PDPs. PDPs constitute an important family of urban logistics, where goods are transported from various origins to corresponding destinations. The PDP can be categorized into three primary groups based on the nature of the origin-destination pairs: one-to-one problems, where requests start and end at one origin and one destination; one-to-many problems, involving the transport of items from one origin to multiple destinations; and many-to-many problems, where multiple entities (or commodities) are transported between numerous origins and destinations \citep{daganzo2005logistics}. The one-to-one distribution problem primarily aims to optimize the one-dimensional dispatching headways along the timeline in response to dynamic demand, which is traditionally addressed through dynamic programming \citep{newell1971dispatching}. The one-to-many and many-to-many problems are grounded in classic routing problems, such as the Traveling Salesman Problem (TSP) and the Vehicle Routing Problem (VRP), and they vary based on the number of vehicles deployed. 

To estimate the minimum total traveling distance for one-to-many TSP, CA models suggest that the optimal distance can be approximated as $d_{TSP} \approx k_{TSP}\sqrt{AN}$, where $A$ represents the area of service region, and $N$ denotes the number of destinations to visit within this area. Through a swath heuristic, which splits the service area into several strips and then determines the optimal swath width, \cite{daganzo1984length} estimated the constant factor $k_{TSP}$ is 0.9 for Euclidean metric and 1.15 for $L_1$ metric if the service zone is not too narrow. Different from the swath-split strategy, \cite{del1999heuristic} constructs the suboptimal tour via a ring-radial separation. Their numerical experiments demonstrate that the average Euclidean length of the tours produced by this approach scales with $\sqrt{AN}$, showing performance essentially on par with Daganzo's swath heuristic. Building on the TSP approximation, the multiple-vehicle delivery problem, or capacitated vehicle routing problem, typically assumes each vehicle can only deliver to a limited number of points. The expected all-vehicle travel distance is expressed similarly to the TSP as $d_{VRP} \approx 2r(s/S)N+k_{VRP}\sqrt{AN}$, where $r$ represents the distance from the depot to the service area's center, $s$ and $S$ are demand expectation and vehicle capacity, respectively \citep{daganzo1984distance}. This formula, incorporating a 'line-haul distance' for initial travel from the depot to the customers and a 'local distance' for detours within each region, is also extended to a ring-radial network and other rectangular zones with various orientations \citep{daganzo1986design, newell1986design}.

The many-to-many distribution problem, usually referred to as the multi-commodity problem, assumes each demand has a specific origin and destination. The many-to-many distribution system can be divided into two categories: VRPPD and dial-a-ride problem (DARP), where the former considers the transportation of goods, and the latter focuses on passenger transportation. The key difference is that DARP usually incorporates an extra hard or soft constraint to specifically address and minimize customers' inconvenience \citep{toth2002vehicle}. \cite{daganzo1978approximate} initially formulated a heuristic strategy designed to estimate customers' waiting and riding times, which can be represented as functions of the number of vehicles, request density, service area, and vehicle traveling speed. The optimal approach advocates for alternating pickups and deliveries, consistently opting for the nearest pick-up or delivery point. The effectiveness of this approach was illustrated through a comparison against two alternative strategies: (i) the bus routes to the nearest feasible point, and (ii) the bus collects a fixed number of passengers before proceeding with deliveries. To estimate the total service time, several simulation models \citep{wilson1969simulation, wilson1976advanced} and empirical models \citep{flushberg1976descriptive} have also been developed for the many-to-many DARP service. \cite{stein1978asymptotic} proved that the many-to-many TSP with pickup and delivery constraints (TSPPD), adheres to a similar formula as the traditional TSP with $d_{TSPPD} \approx k_{TSPPD}\sqrt{AN}$, where $\sqrt{2}k_{TSP} \leq k_{TSPPD} \leq 2k_{TSP}$. Through a probabilistic analysis, Stein split the optimal tour into several segments and demonstrated that the optimal $k_{TSPPD}$ asymptotically converges to $4\sqrt{2}/3$. In this study, we adopt the methodology from \cite{daganzo1978approximate} by incorporating additional practical considerations specific to meal delivery to estimate deliverers' VMT and customers' waiting time under the VRPPD strategy. These estimations are then used to benchmark the performance of a distribution system that implements transshipment.

\subsection{Distribution system with transshipment}

In one-to-many or many-to-many distribution, transshipment can be viewed as a terminal sorting, storing, and transferring items from one vehicle to another. Incorporating transshipment can significantly reduce the total miles traveled by vehicles. But at the same time, it may introduce new handling and holding costs at the terminal, attributable to the increased batch sizes. Designing the distribution system with transshipment typically involves determining the number of terminals to be utilized, their strategic locations, the planning of routes and schedules for various types of vehicles, and the assignment of customers to specific terminals and routes \citep{daganzo2005logistics}. \cite{daganzo1986configuration} first employed a CA method to find near-optimal solutions for the design guidelines in a one-to-many distribution system, such as the location of transshipping points, their areas of influence, vehicle frequency, and fleet size. \cite{smilowitz2007continuum} then advanced the design strategy for an integrated package distribution system considering service level, i.e., overnight and longer deadlines. They separated the complex system into a series of subproblems to determine the densities of various terminal levels (e.g., consolidation terminal, breakbulk terminal, airport) and their corresponding dispatching headways through a CA method. Extensive studies have employed the CA method to design distribution systems or identify the optimal location for transshipment across various application contexts, such as delivery consolidation in freight transportation \citep{campbell2013continuous, xie2015optimal, ghaffarinasab2018continuous}, inventory management within supply chain \citep{naseraldin2011location, tsao2012continuous}, and agricultural collection during harvest season \citep{wiles2001optimal}. 

To evaluate the benefit of the transshipment design, \cite{daganzo1987break} first delved into the `break-bulk' role of the transshipment in many-to-many logistic networks using CA. The investigation highlighted the cost-saving benefits of allowing one and then two transshipment stops, attributing these savings to the reduction in line-haul distance traveled. They offset the additional costs associated with transshipment, including delivery delays and handling costs at the terminal. In the modern logistics network with multiple suppliers and customers, numerous studies have employed CA to evaluate the effects of delivery consolidation. To reduce transportation costs in urban delivery, companies collaborate to consolidate shipments at a central terminal rather than independently dispatching goods from the distribution centers to customers. \cite{kawamura2007evaluation} analyzed the economic efficiency of cooperative delivery consolidation based on a logistics cost analysis and concluded that, without accounting for societal benefits such as reduced illegal parking and congestion, the appeal of consolidation to the U.S. industry is marginal. \cite{chen2012comparison} extended this work by considering a different network setting for small businesses and using an additional dataset for parameter estimation. Their study theoretically demonstrated that delivery consolidation becomes cost-effective when the number of customers and suppliers is large, and the terminal operation cost is low. \cite{lin2016sustainability} then expanded their CA model by incorporating energy consumption and PM2.5 emissions as additional costs, along with policy factors like commercial vehicle size restriction in urban areas. Their findings indicate that when there is an economy of scale or high customer density, the consolidation design yields both logistics and environmental benefits. Similarly, \cite{pahwa2022cost} developed a multi-echelon last-mile distribution model and suggested that consolidation strategies are well-suited for areas with dense demand. In contrast, for regions with lower population density, outsourcing solutions such as crowdsourced delivery and customer self-collection emerge as more cost-effective alternatives. 

In urban logistics, existing research on transshipment has mostly focused on the design of consolidation centers for logistic companies, with an emphasis on the operational costs, such as handling and fixed pipeline inventory expenses associated with additional transshipment. However, these costs are relatively minor when applied to the temporary storage of meal packages awaiting transshipment. The effectiveness of transshipment strategies for meal delivery services, as well as the optimal system design, has not yet been thoroughly investigated and remains unclear. We thus propose to examine the cost-effectiveness of transshipment specifically for meal delivery from the perspectives of both deliverers and customers. In this study, we develop stylized models for meal delivery operations, both with and without transshipment, to compare the efficiency of the two operational models.

\section{Integrated Meal Delivery System with Transshipment}\label{sec:microhub}

This section introduces our proposed system design for transshipment and key performance metrics related to different stakeholders. It begins with an overview of the system setup for the transshipment operations, followed by a breakdown of the lifecycle for meal orders from creation to package fulfillment, from their creation to final package delivery. Given the system inputs, continuous approximations are applied to estimate the average vehicle miles traveled by deliverers, whilst a queuing network is conceptualized to analyze the time that meal packages spend at various delivery stages.

\subsection{Meal delivery with transshipment}

Consider a meal delivery service offered in a region of area $A$, with an average order arrival flux $\lambda$ (i.e., order arrivals per hour per square mile), a total of $m$ deliverers. For the transshipment operations illustrated in Figure \ref{fig_intro_d}, the region is divided into $K$ pie-shaped sub-areas, and each deliverer is assigned to service within a specific sub-area. Each sub-area $k\in K$ is characterized by a coverage area $A_k$, an order arrival rate $\lambda_k(=\lambda A_k)$, and a dedicated fleet of $m_k$ deliverers. Any regional partition defined by $\{A_k, m_k\}_{k\in K}$ must satisfy the constraints $\sum_k m_k = m$ and $\sum_k A_k =  A$. Given a partition of sub-areas, individual delivery orders can be categorized as either intra-zonal or inter-zonal. Intra-zonal orders have both their pickup and drop-off locations within the same sub-area, whereas inter-zonal orders involve pickup and drop-off locations spanning different sub-areas.

Deliverers operate cyclically and exclusively within their designated sub-areas to perform pickup and drop-off tasks. In each sub-area $k$, a deliverer is dispatched once there are $n_k$ points accumulated to visit. These $n_k$ points comprise both the meals waiting for drop-off at the microhub and the meals awaiting pickup at restaurants. Upon dispatch, the deliverer completes a round trip that begins and ends at the microhub, touring through the $n_k$ pending pickup and drop-off points. Regardless of whether the orders are intra-zonal or inter-zonal, all meal packages are first picked up from their originating restaurants and transported to the microhub during the current delivery cycle. At the microhub, they are sorted according to their destinations and await drop-off. Then, in a subsequent cycle, one deliverer assigned to the sub-area containing an order's drop-off point retrieves the package from the microhub and delivers it to its final destination. As a result, the distribution tasks in each delivery cycle within a sub-area closely resemble a one-to-many-to-one (1-M-1) pickup and delivery problem, involving two distinct sets of commodities: some are transported from the microhub to customers, while others are picked up at restaurants and returned to the microhub \citep{toth2002vehicle}.

\subsection{Lifetime of meal delivery orders}\label{sec:order_lifetime}

\begin{figure}[!t]
	\centering
	\includegraphics[width=0.9\textwidth]{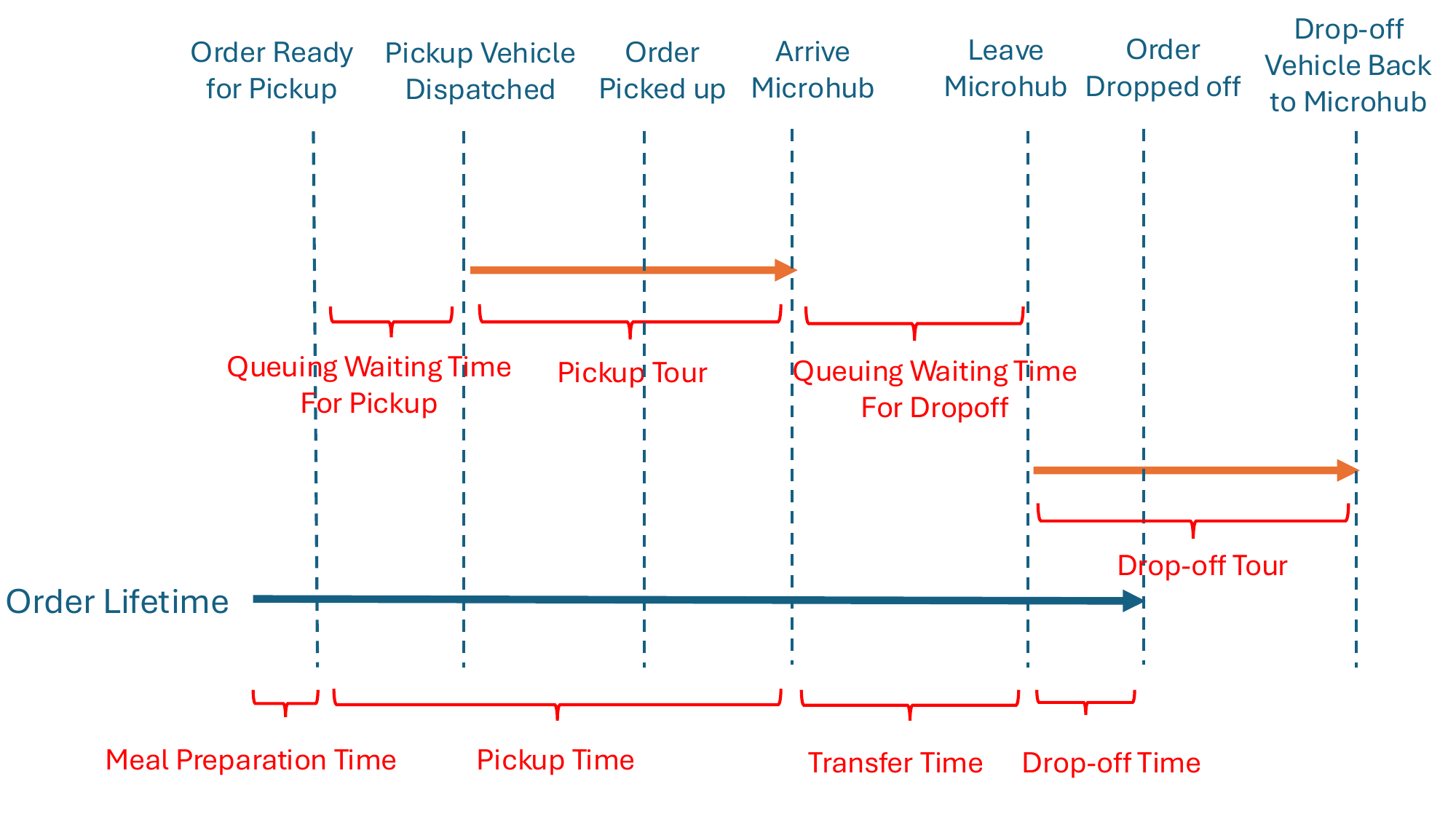}
	\caption{Lifetime of individual meal delivery orders}
	\label{fig_order_lifetime}
\end{figure}

Figure \ref{fig_order_lifetime} outlines the lifecycle of a delivery order for meal packages under the proposed transshipment strategies. Each delivery order progresses through four sequential stages, defined as follows:
\begin{enumerate}
\item \emph{Meal Preparation}: Interval from when an order is placed until it is ready for pickup;
\item \emph{Pickup}: Interval from the meal package being ready to its arrival at the microhub;
\item \emph{Transfer}: Interval during which the package is held at the microhub for the next delivery tour;
\item \emph{Drop-off}: Interval for transporting the package from the microhub to the customer.
\end{enumerate}

\noindent Of these, the meal preparation stage is unaffected by the introduction of transshipment and will therefore be left out in the following analysis. The pickup stage can be further broken down into three phases: the period from when the package is ready for pickup until a deliverer is dispatched, the interval between the dispatch of the deliverer and the actual pickup of the order, and the duration required for the package to be transported from the restaurant to the microhub. The combined duration of the latter two phases indeed constitutes the total time for a complete delivery tour executed by a deliverer. Additionally, it should be noted that when a meal is ready for pickup, it does not guarantee the immediate dispatch of a deliverer. Two conditions must be met within each sub-area: first, $n_k$ pending pickup and drop-off tasks must accumulate in sub-area $k$; second, after $n_k$ tasks have been accumulated, the batch must wait until a deliverer assigned to sub-area $k$ becomes available to handle the service. Due to the presence of inter-zonal orders, the operations across different sub-areas are closely interdependent. To capture these cross-zonal interactions, a queuing network with $K$ dependent queues is proposed to pin down the waiting time for pickup as well as the transfer time, which can also be interpreted as the waiting time for drop-off.

\subsection{Queuing network model}\label{sec:queuing_network}

The $K$ sub-areas collectively form a Jackson queuing network in Figure \ref{fig_queuing_network}, where each sub-area operates as a separate queue with its own servers. The following assumptions are made for this network of $K$ single-station queues:

\begin{figure}[ht!]
	\centering
	\includegraphics[width=0.9\textwidth]{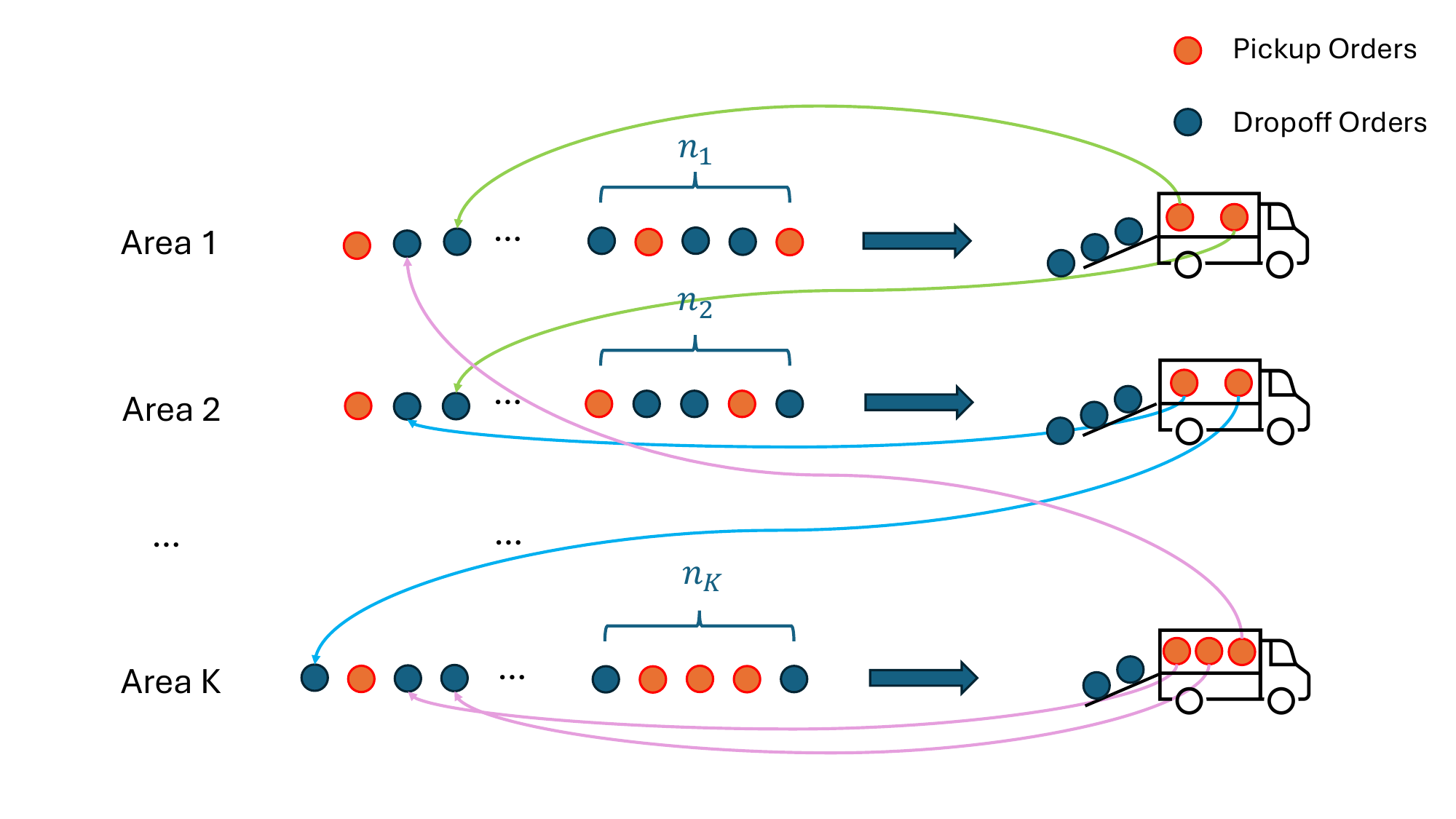}
	\caption{Cross-zonal meal delivery process as a queuing network}
	\label{fig_queuing_network}
\end{figure}

\begin{itemize}
    \item New orders arrive at sub-area $k$ following a Poisson process with an arrival rate of $\lambda_k$. The pickup point for each order is randomly located within sub-area $k$.
    \item A package picked up from sub-area $i$, after being transferred to the microhub, will join the queue at sub-area $k$ for drop-off with probability $p_{i,k}$. The drop-off points of orders are also random within each sub-area $k$. Once a package is dropped off, it exits the queuing network completely.
    \item The $k$-th sub-area is served by $m_k$ delivery vehicles. When the number of orders in the queue reaches $n_k$ or more and a deliverer becomes idle, a batch of exactly $n_k$ orders is taken for service.
    \item The service time $S_k$ for meal packages at sub-area $k$ are independent and identically distributed, with a mean $\E[S_k]$ and a variance $\Var[S_k]$.
    \item The microhub has sufficient capacity to accommodate all incoming packages, and the queue at each sub-area has unlimited waiting spots.
\end{itemize}

The arrival of meal orders in the queue for each sub-area can be viewed as a Poisson process with an arrival rate $\delta_k$, which sums the arrival rates of newly generated orders in sub-area $k$ for pickup and the collected orders from all sub-areas that need to be dropped off to sub-area $k$, i.e.,
\begin{align}\label{eq_delta}
    \delta_k = \lambda_k + \sum_{i=1}^K \lambda_i p_{i,k}.
\end{align}
The queued orders are serviced in batches of size $n_k$. For the $i$-th order arriving before reaching the $n_k$ threshold in a single batch, its waiting time is given by $(n_k-i)/\delta_k$. Consequently, the average waiting time $W_k^a$ of an arbitrary pickup or drop-off order in sub-area $k$, from its arrival until the accumulation of $n_k$ orders for a batch, is calculated as follows,
\begin{align}\label{eq_W_a}
    W_k^a &= \frac{1}{n_k} \sum_{i=1}^{n_k} \frac{n_k-i}{\delta_k} = \frac{n_k-1}{2\delta_k}.
\end{align}

After accumulating a sufficient number of orders for a batch, the batched orders may also experience additional waiting time if a deliverer is not immediately available. For each sub-area, this process of forming and delivering batch orders by deliverers can be modeled as a $G/G/m$ queue. The inter-arrival time of consecutive order batches follows an Erlang distribution with the shape parameter $n_k$ and the rate parameter $\delta_k$. The service rate by a deliverer is determined by the touring time required to complete a Traveling Salesman Problem (TSP) for fulfilling the batched orders. The number of servers in each queue equals the number of deliverers $m_k$. Thus, the additional queuing time for batched orders $W_k^q$ can be estimated with a widely used approximation for a $G/G/m$ queue by \cite{larson1981urban}, presented as follows,
\begin{align}\label{eq_W_q}
    W_k^q \approx \frac{\left[ \frac{n_k}{\delta_k^2}+\frac{1}{m_k}\Var[S_k] \right] \cdot \frac{\delta_k}{n_k}}{2 \left(1-\frac{\delta_k \cdot \E[S_k]}{n_k \cdot m_k}\right)}.
\end{align}

\subsection{Service time approximations}\label{sec:TSP_approx}

The service process for a single queue follows a TSP tour undertaken by deliverers. \cite{daganzo1984length} showed that, in a regularly shaped area of size $A$ with $N$ randomly distributed points to visit, the expected TSP tour length is proportional to $\sqrt{AN}$. To adapt this formulation for the pie-shaped sub-areas, we introduce an extra distance $R'$ to account for radial forward and return travel between the microhub and the farthest points within each pie-shaped region. Proposition \ref{proposition_1} derives the expected value of $R'$.
\begin{prop}\label{proposition_1}
    For a circular sector of radius $R$ containing $N$ nodes, the expected farthest distance, $\E[R']$, is given by $\E[R'] = \frac{2N}{2N+1}R$, where $\E[R']$ is independent of the sector's angle.
\end{prop}
\noindent This formula suggests that when there are fewer points to visit, the radial traveling distance often remains short, as deliverers likely do not need to venture far from the microhub. However, as the number of points increases, deliverers tend to travel closer to the outer perimeter of the pie-shaped area. By incorporating this radial distance, the expected service time in sub-area $k$ under the modified TSP approximation is expressed as follows:
\begin{align}
    \E[S_k] &= \frac{1}{v}\cdot (a \cdot \sqrt{A_k n_k} + b \cdot \E[R'_k]) \label{eq_E_Dk}, \\
    % \E[S_k] &= \E[\frac{D_k}{v}] \label{eq_E_Sk}, \\
    \E[R'_k] &= \frac{2n_k}{2n_k+1}R_k \label{eq_R},
\end{align}
where $v$ denotes the speed of delivery vehicles, and $R_k$ is the radius of sub-area $k$, which approximates the average radial distance from the vertex to the perimeter border of the pie-shaped area. It is worth noting that the origin-destination pairs of individual meal delivery orders may be correlated according to certain distance profiles, implying that pickup and drop-off locations are not entirely random within the service region. However, in the transshipment system, each order’s pickup and drop-off points are separated into two distinct delivery tours. As a result, all pickup and drop-off points within a single tour correspond to different orders. Consequently, it is reasonable to assume that the access points visited by individual deliverers within each sub-area are randomly distributed in space.

For the variance of TSP touring distances $\Var[S_k]$, it is intuitive that its value is linearly proportional to the size $A_k$ of the sub-area. Regarding its relationship with the number of points $n_k$, we observe that when $n_k$ is small, the spatial dispersion of points can be highly uncertain, with points either clustering together or spreading across the entire region. Conversely, as $n_k$ increases, the points tend to distribute more uniformly throughout the region, leading to a more stable variance. Thus, $\Var[S_k]$ could be inversely proportional to $n_k$ and converge to a certain positive value as $n_k$ approaches infinity. Based on the above deduction, we propose the following hypothetical form for describing the service time variance under the TSP:
\begin{align}
    \Var[S_k] &= \frac{A_k}{v^2} \cdot \left( \frac{\alpha}{n_k} + \beta \right) \label{eq_Var_Dk}.
    % \Var[S_k] &= \Var[\frac{D_k}{v}] \label{eq_Var_Sk}.
\end{align}
The hyperparameters $a, b$ and $\alpha, \beta$ in Eqs. \eqref{eq_E_Dk} and \eqref{eq_Var_Dk} will be calibrated later through simulation, where the effectiveness of the proposed approximations will also be demonstrated.

\subsection{System performance metrics}

Given a regional partition $\{A_k, m_k\}_{k\in K}$ and the order batching sizes $\{n_k\}_{k\in K}$, the queuing network model predicts the stationary states of the meal delivery system under transshipment and subsequently yields the two key system performance metrics.

\textit{Deliverers' vehicle miles traveled.}\quad For each sub-area, when $n_k$ orders are collected, one deliverer is dispatched to visit these $n_k$ points, while other idle deliverers, if any, wait near the microhub waiting for the next delivery. Thus, the expected VMT per hour by all deliverers in sub-area $k$, denoted as $Q_k$, is given by the product of the number of delivery cycles performed per hour and the expected traveling distance for visiting $n_k$ points. Specifically, 
\begin{align}\label{eq_miles_traveled}
    Q_k &= \frac{\delta_k}{n_k} \cdot \E[S_k]\cdot v.
\end{align}

\textit{Customers' waiting time.}\quad As previously noted, each customer's waiting time can be divided into four components, with the meal preparation time excluded from this analysis. Combining the rest three components, the expected total waiting time for meal orders delivered from sub-area $i$ to $k$, denoted as $W_{i,k}$, can be expressed as:
\begin{align}\label{eq_waiting_time}
    W_{i,k} &= W_i^a + W_i^q + \E[S_i] + W_k^a + W_k^q + \frac{\Var[S_k] + \E^2[S_k]}{2 \E[S_k]},
\end{align}
where $W_i^a + W_i^q + \E[S_i]$ accounts for the pickup time, covering the period from order generation to arrival at the microhub; $W_k^a + W_k^q$ represents the transfer time at the microhub; and the last term in Eq. \eqref{eq_waiting_time} stands for the drop-off time. Specifically, the time for package drop-off is derived based on the occurrence of a random event along a TSP tour, where the touring time follows a given distribution with a known mean and variance \citep{larson1981urban}.

\section{A VRPPD Benchmark without Transshipment}\label{sec:VRPPD}

For the benchmark operations without transshipment, individual deliverers continuously traverse the service area, picking up meal packages and delivering them directly to their destinations. In this section, we adopt a heuristic for the VRPPD, originally proposed by \cite{daganzo1978approximate}, to estimate system performance without transshipment. The heuristic is further enhanced by incorporating practical considerations specific to the delivery distances of individual meal orders. 

\subsection{A heuristic for meal delivery without transshipment}

The heuristic devised by \cite{daganzo1978approximate} directs deliverers to continuously route to the nearest point pending a visit, which can be either the origin of a pickup order or the destination of a drop-off order already on the vehicle. At equilibrium, a deliverer is equally likely to visit an origin/pickup point or a destination/drop-off point. By assuming that both origins and destinations are independently and randomly distributed across the service area, \cite{daganzo1978approximate} showed that the trip chains continue with equal probabilities, and the number of packages on each delivery vehicle simply equals the number of packages waiting for pickup. Unlike the setting in \cite{daganzo1978approximate}, distance-based delivery fees and platform-imposed limitations on ordering distances create strong bundling between origins and destinations of meal orders, shaping delivery distances into specific profiles. The empirical data from Meituan calibrates the distance profile as a Rayleigh distribution, with the goodness-of-fit detailed in Appendix \ref{sec:appd_Rayleigh}.

\begin{figure}[b!]
	\centering
	\begin{subfigure}[b]{.5\textwidth}
		\centering
		\includegraphics[width=0.9\linewidth]{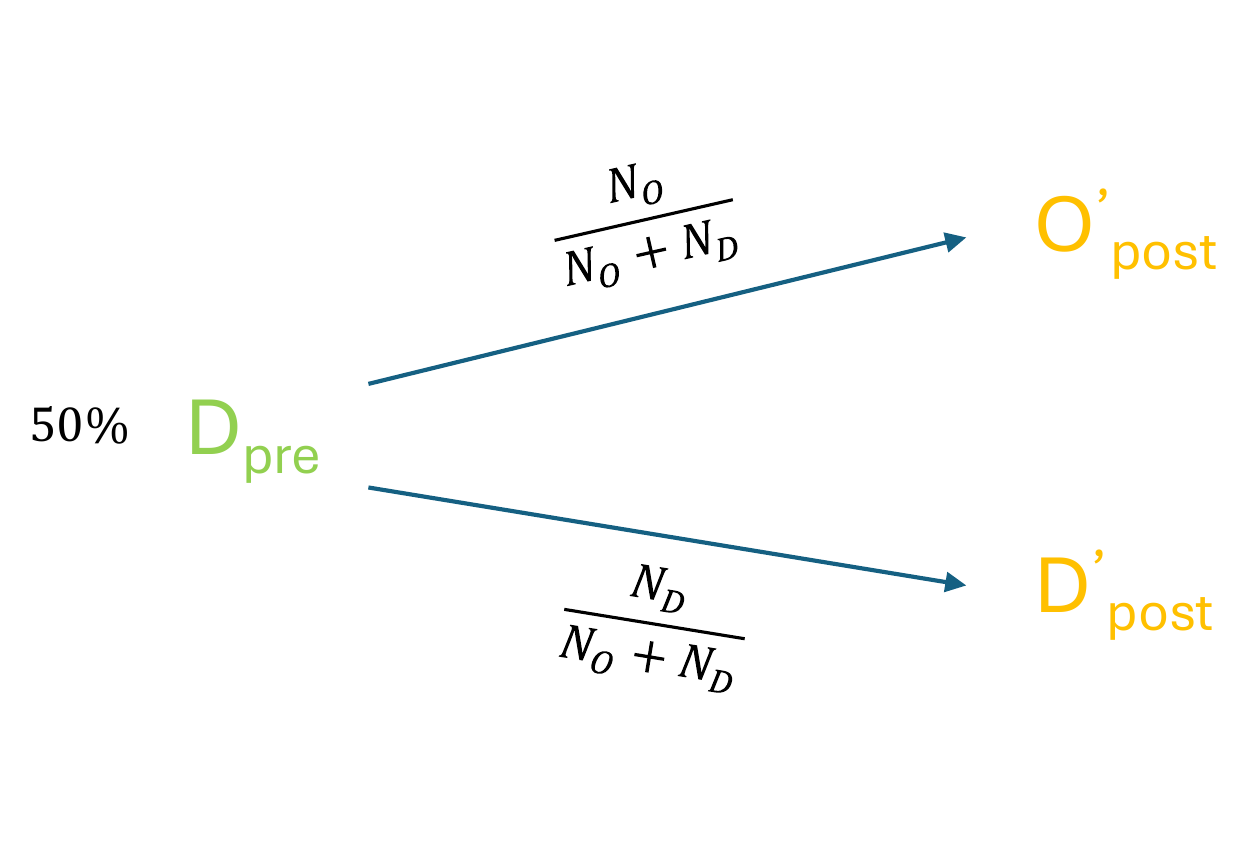}
		\caption{}
        \label{fig_VRPPD_D}
	\end{subfigure}%
	\begin{subfigure}[b]{.5\textwidth}
		\centering
		\includegraphics[width=0.9\linewidth]{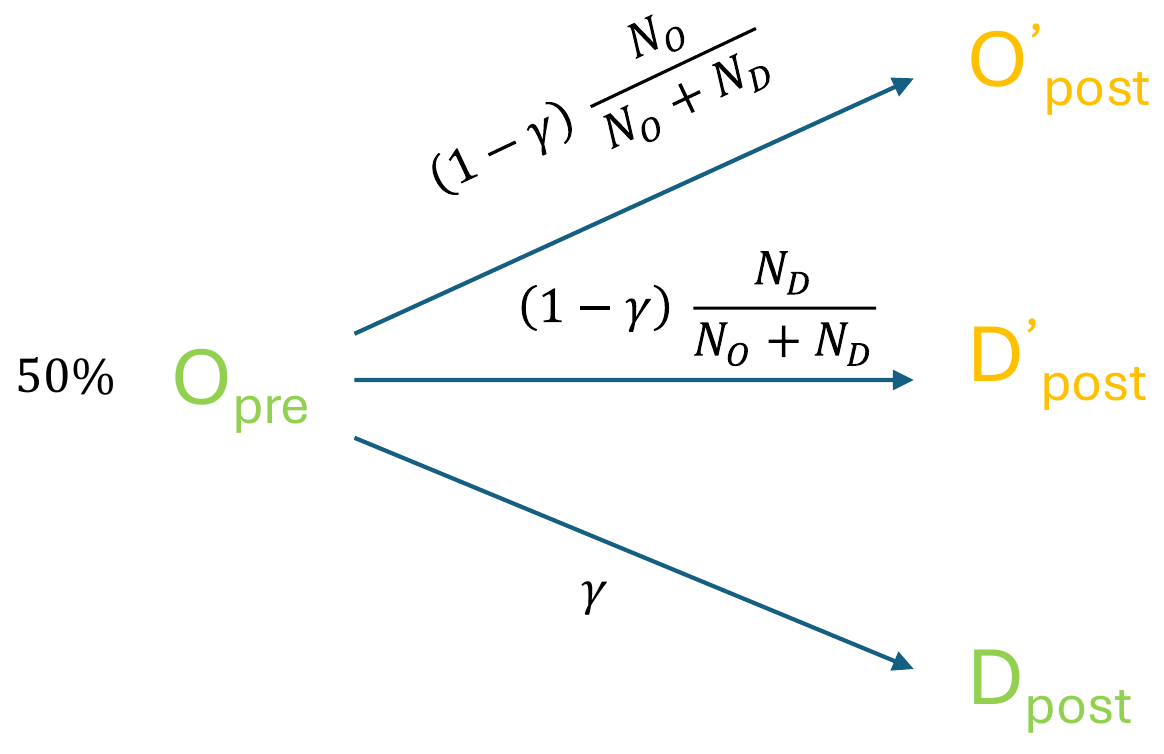}
		\caption{}
        \label{fig_VRPPD_O}
	\end{subfigure}
        \caption{Scenario enumeration of state transitions between the types of preceding and succeeding points visited consecutively by individual deliverers.}
	% \caption{Demonstration of the customized VRPPD strategy, illustrating the current visit to (a) a destination point and (b) an origin point.}
	\label{fig_VRPPD}
\end{figure}

Accounting for the distance bundling between origin-destination pairs of meal delivery results in more sophisticated state transition probabilities for deliverers' trip chains, as illustrated in Figure \ref{fig_VRPPD}. Let $N_O$ denote the total number of packages waiting for pickup and $N_D$ the number of packages en-route on a delivery vehicle. Then, immediately following a visit to a drop-off point $D_{pre}$, shown in Figure \ref{fig_VRPPD_D}, the deliverer has a probability of $\frac{N_O}{N_O + N_D}$ to visit a pickup point $O'_{post}$ and $\frac{N_D}{N_O + N_D}$ to visit another drop-off point $D'_{post}$. These probabilities arise from the assumption that all meal orders are generated i.i.d. across the region according to a spatial Poisson process. On the other hand, following a visit to a pickup point $O_{pre}$, Figure \ref{fig_VRPPD_O} shows that the deliverer's next visit falls into either of the two cases: (1) the corresponding drop-off point $D_{post}$, with probability $\gamma$, or (2) any other access point within the service region, with probability $1 - \gamma$. The latter case can be further decomposed into visits to either another pickup point $O'_{post}$ or a drop-off point $D'_{post}$ that is different from the destination of the order just picked up. The probability $\gamma$ represents the likelihood that, after visiting a pickup point, the distance to its corresponding drop-off point is shorter than to any other point to visit. This probability can be derived from Proposition \ref{proposition_2}.

\begin{prop}\label{proposition_2}
    Let $X$ be a random variable representing the distance to the nearest point that adheres to a spatial Poisson process with density $\frac{N_O+N_D}{\pi R^2}$. Let $Y$ be a random variable following a Rayleigh distribution with a scale parameter $\sigma>0$. Then, the probability that $X$ greater than $Y$ is given by:
    \begin{align}\label{eq_VRPPD_p}
        \gamma = \probP(X>Y) = \frac{1}{1+\frac{2(N_O+N_D)\sigma^2}{R^2}}.
    \end{align}
\end{prop}
\noindent As suggested by Eq. \eqref{eq_VRPPD_p}, if the delivery orders are shorter, with a smaller $\sigma$, the probability of directly visiting the drop-off point immediately after picking up an order without detours increases. Conversely, with denser visiting points in space, the likelihood of immediately dropping off an order following its pickup decreases.

The state transitions spelled out in Figure \ref{fig_VRPPD} lead to the following relationship based on the Law of total probability:
\begin{align*}
    \probP(O_{post}) &= \probP(D_{pre}) \cdot \probP(O_{post} | D_{pre}) + \probP(O_{pre}) \cdot \probP(O_{post} | O_{pre}), \\
    \Leftrightarrow\quad \probP(O_{post}) &= \probP(D_{pre}) \cdot \frac{N_O}{N_O+N_D} + \probP(O_{pre})  \cdot (1-\gamma) \frac{N_O}{N_O+N_D} .
\end{align*}
% and
% \begin{align*}
%     \probP(D_{post}) &= \probP(D_{pre}) \cdot \probP(D_{post} | D_{pre}) + \probP(O_{pre}) \cdot \probP(D_{post} | O_{pre}) \\
%     \Leftrightarrow\quad 50\% &= 50\% \cdot \frac{N_D}{N_O+N_D} + 50\% \cdot \left( \gamma + (1-\gamma) \frac{N_D}{N_O+N_D} \right) 
% \end{align*}
Moreover, the equilibrium condition implies that $\probP(O_{pre}) = \probP(D_{pre}) = \probP(O_{post}) = \probP(D_{post}) = 50\%$, which then gives rise to the following equation:
\begin{align}\label{eq_VRPPD_pN}
    \gamma = \frac{N_O - N_D}{N_O}.
\end{align}

\subsection{Service rate estimations}

The service rate $\mu$ provided by individual deliverers is given by:
\begin{align}
    \mu &= \left(\frac{2 \E[d]}{v}\right)^{-1}, \label{eq_VRPPD_mu}
\end{align}
where $d$ is the distance from an arbitrary visit point to its nearest point. To derive this, we analyze the two scenarios separately to determine the distances $d_D$ and $d_O$, which respectively denote the distance to the nearest point after visiting a drop-off and a pickup point, given $N_O + N_D$ points to visit within the service area. 

First, with the intensity of the spatial Poisson process given by $\delta = \frac{N_O + N_D}{\pi R^2}$, the cumulative distribution function (CDF) and probability density function (PDF) of $d_D$ are $F_{d_D}(r) = 1-e^{-\delta \cdot \pi r^2}$ and $f_{d_D}(r) = 2 \pi \delta r \cdot e^{-\delta \cdot \pi r^2}$. The expected value of $d_D$ is thus,
\begin{align}
    \E[d_D] & = \int_0^{+\infty} r \cdot 2 \pi \delta r \cdot e^{-\delta \cdot \pi r^2} \, \d r \nonumber\\
    & = \int_0^{+\infty} e^{-\delta \cdot \pi r^2} \, \d r \nonumber\\
    & = \frac{1}{2} \sqrt{\frac{\pi}{\delta\cdot\pi}} \nonumber\\
    & = \frac{\sqrt{\pi}R}{2\sqrt{N_O + N_D}}.
\end{align}
Specifically, the first equality follows from the definition, the second applies the integration by parts, the third uses the Gaussian integral formula, and the final equality substitutes the corresponding expression for $\delta$.
% Let $u=r$ and $v=\int 2 \pi \delta r \cdot e^{-\delta \cdot \pi r^2} dr = -e^{-\delta \cdot \pi r^2}$. Applying the formula for integration by parts,
% \begin{align*}
%     \E[d_D] &= \left[ r \cdot \left(-e^{-\delta \cdot \pi r^2} \right) \right]_0^{+\infty} - \int_0^{+\infty} -e^{-\delta \cdot \pi r^2} \, dr \\
%     &= \int_0^{+\infty} e^{-\delta \cdot \pi r^2} \, dr.
% \end{align*}
% Substitute $s=\delta \pi$, the integral simplifies using the Gaussian integral formula,
% \begin{align*}
%     \int_0^{+\infty} e^{-sr^2} \, dr = \frac{1}{2} \sqrt{\frac{\pi}{s}}.
% \end{align*}
% Thus, we find:
% \begin{align}
%     \E[d_D] &= \frac{1}{2\sqrt{\delta}} = \frac{\sqrt{\pi}R}{2\sqrt{N_O + N_D}}.
% \end{align}

Second, immediately after visiting a pickup point, the minimum distance to the next point is determined by the closest point between its corresponding drop-off point and other existing visit points, as derived in Proposition \ref{proposition_3}. 
\begin{prop}\label{proposition_3}
    As a continuation of Proposition \ref{proposition_2}, the expected value of the minimum of $X$ and $Y$ is given by $\E[\min(X,Y)] = \frac{1}{2} \sqrt{\frac{\pi}{\frac{N_O+N_D}{R^2}+\frac{1}{2\sigma^2}}}$.
\end{prop}

Combining the two cases, the average distance $d$ from an arbitrary visit point to its nearest point is obtained as follows,
\begin{align}\label{eq_VRPPD_d}
    \E[d] &= \probP(D_{pre})\cdot \E[d_D] + \probP(O_{pre})\cdot \E[d_O] = \frac{\sqrt{\pi}R}{4\sqrt{N_O+N_D}} + \frac{1}{4} \sqrt{\frac{\pi}{\frac{N_O+N_D}{R^2}+\frac{1}{2\sigma^2}}}.
\end{align}

\subsection{System performance metrics}
When the demand for meal delivery is low, it is unnecessary to have all deliverers actively in service. Therefore, we define an additional variable $m'$ as the total number of deliverers in service, which should not exceed the total number of available deliverers $m$. Additionally, system equilibrium requires that the arrival rate ($\lambda A$) equals the total service rate ($\mu m'$), i.e.,
\begin{align} \label{eq_VRPPD_equilibrium}
    \lambda A &= \mu m'.
\end{align}
Combining all the components above, given the demand flux $\lambda$ and the number of deliverers in service $m'$, the equilibrium state of the VRPPD system without transshipment can be determined by solving the system of equations \eqref{eq_VRPPD_p}-\eqref{eq_VRPPD_equilibrium}.

\emph{Deliverers' vehicle miles traveled}.\quad Similar to Eq. \eqref{eq_miles_traveled}, the total VMT by deliverers is given by the product of the number of points visited by a deliverer per hour and the expected distance between two adjacent visits:
\begin{align}\label{eq_VRPPD_Q}
    Q_{VRPPD} &= 2\lambda A \cdot \E[d] = m' \cdot v.
\end{align}

\emph{Customers' waiting time}.\quad The average waiting time for customers, including both the waiting time for package pickup and the riding time on a vehicle, is given by:
\begin{align}\label{eq_VRPPD_W}
    W_{VRPPD} &= \frac{N_O}{\lambda A} + \frac{N_D}{\mu}.
\end{align}

% Therefore, the VRPPD optimization problem below is to determine the optimal $m'$, along with the corresponding values of $N_O$ and $N_D$,
% \begin{subequations}\label{opt_VRPPD}
% 	\begin{align}
% 		& \min_{\{m', N_v, N_w,\}} \pi_Q \cdot Q_{VRPPD} + \pi_W \cdot \lambda A \cdot W_{VRPPD} \\
% 		\text{s.t.}\ \ \ & (\ref{eq_VRPPD_p}) - (\ref{eq_VRPPD_Q}) \nonumber \\
%             & N_v, N_w, m' \geq 0
% 	\end{align}
% \end{subequations}

\section{Numerical Analysis}\label{sec:numerical}

We then evaluate the potential of the proposed meal delivery operations with transshipment by benchmarking them against the traditional pickup-and-delivery model without transshipment. First, we set the stage for numerical experiments by formulating two optimal system designs to establish the performance envelope under different market conditions. Next, we run simulations to validate the accuracy of our analytical model in predicting average customer waiting time and VMT by deliverers. A comprehensive performance evaluation is then conducted to compare meal delivery systems with and without implementing transshipment. Finally, we introduce a dataset from the Meituan meal delivery platform to compare the performance of these strategies in a real-world scenario.

\subsection{Experimental setup}

To ensure a fair head-to-head comparison of operations using the proposed CA models, this evaluation considers the context of an isotropic circular area of radius $R$ with an average order arrival flux $\lambda$ and a total of $m$ deliverers. For each market condition, we optimize the system design for both transshipment and non-transshipment operations and compare the resulting system performance.

For transshipment operations, the experiments in this section focus on symmetric partitioning strategies, where the area is evenly divided into $K$ sub-areas, with deliverers uniformly distributed across them. For brevity, we retain the notations of variables previously defined for individual sub-areas with the subscript $k$. However, since all sub-areas are now identical, these notations will now refer to a representative sub-area rather than a specific one. The optimal design under symmetric partitioning reduces to the following optimization program, which determines the optimal number of partition $K$ and the batch size $n_k$ per delivery cycle in sub-areas:
\begin{subequations}\label{opt_form_simplified}
	\begin{align}
		& \min_{K, n_k \geq 0} \ K \cdot \left( \pi_Q \cdot Q_k + \pi_W \cdot \lambda_k W_{k,k} \right) \label{eq_obj}\\
		\text{s.t.}\ \ \ & (\ref{eq_W_a})-(\ref{eq_waiting_time}) \nonumber \\
            & A_k = \frac{A}{K},\ \lambda_k = \lambda A_k,\ \delta_k = 2 \lambda_k,\ m_k = \frac{m}{K} \label{eq_variables_simplified} \\
            & \delta_k \leq \frac{m_k}{\E[S_k]/n_k} \label{cons_steady}
	\end{align}
\end{subequations}
The objective \eqref{eq_obj} represents the generalized variable costs in the meal delivery system, encompassing both the operational costs incurred by deliverers and the time costs experienced by customers, weighted by the coefficients $\pi_Q$ and $\pi_W$. The equalities in \eqref{eq_variables_simplified} specify the allocation of demand and supply within the sub-areas according to the partition. Finally, the condition enforced by \eqref{cons_steady} ensures that the meal delivery system with transshipment remains in steady operations.

For the non-transshipment benchmark modeled as a VRPPD, the optimal system design seeks to determine the number of deliverers actively participating in the service $m'$, along with the total number of packages waiting for pickup $N_O$ and the number of packages en route on a delivery vehicle $N_D$. The corresponding optimization problem is formulated as:
\begin{subequations}\label{opt_VRPPD}
	\begin{align}
		& \min_{m', N_O, N_D \geq 0} \pi_Q \cdot Q_{VRPPD} + \pi_W \cdot \lambda A \cdot W_{VRPPD} \\
		\text{s.t.}\ \ \ & (\ref{eq_VRPPD_p})-(\ref{eq_VRPPD_W}) \nonumber \\
            & m' \leq m \label{cons_active supply}
	\end{align}
\end{subequations}
where constraint \eqref{cons_active supply} ensures that the number of active deliverers in service does not exceed the total available deliverers.

\subsection{Simulation and measurement} \label{sec:simulation_measurement}

Simulations are conducted to demonstrate the effectiveness of proposed approximations for TSP tours---Eqs. \eqref{eq_E_Dk}, \eqref{eq_Var_Dk}---and for customer waiting time---Eqs. \eqref{eq_W_a}, \eqref{eq_W_q}, \eqref{eq_waiting_time}. The simulated results are compared against the predicted approximations across a range of scenarios, varying the values of $\{\lambda, R, m, n, K\}$, to measure the approximation accuracy. Using the trajectory data from Meituan deliverers, the average touring speed for meal delivery is calibrated to $v=4.15$ miles per hour, accounting for delays caused by traffic lights and the time required to access restaurants and customer drop-off points.

\hspace{-1.5em} \textit{Measurement of approximations for deliverers' TSP tours}

\begin{figure}[b!]
	\centering
	\begin{subfigure}[b]{.33\textwidth}
		\centering
		\includegraphics[width=0.95\linewidth]{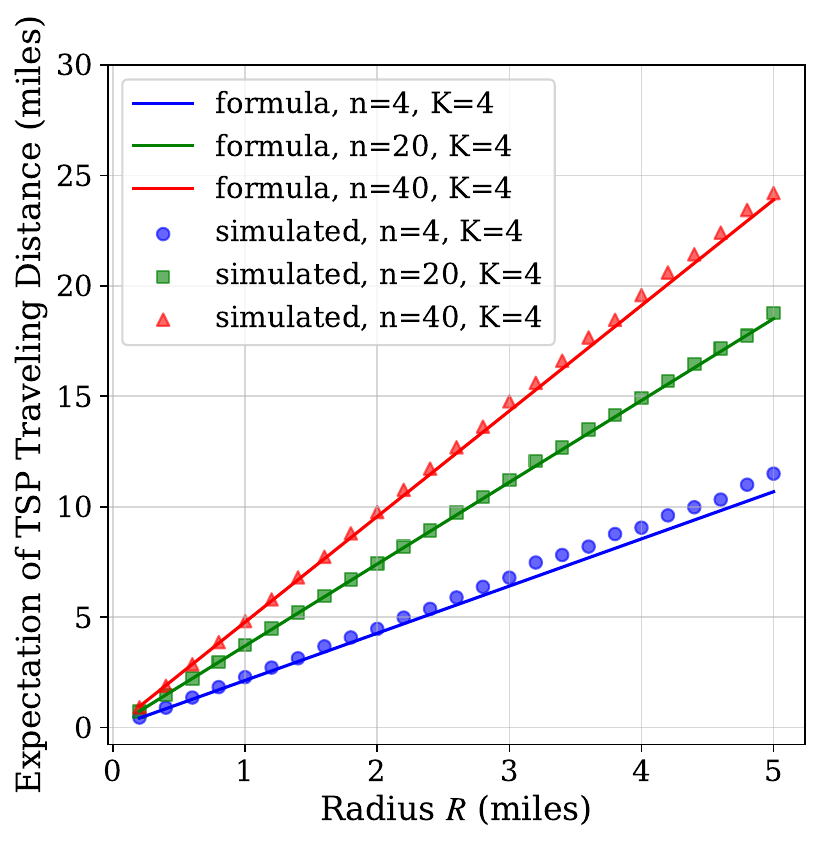}
		\caption{Mean w.r.t. $R$}
	\end{subfigure}%
	\begin{subfigure}[b]{.33\textwidth}
		\centering
		\includegraphics[width=0.95\linewidth]{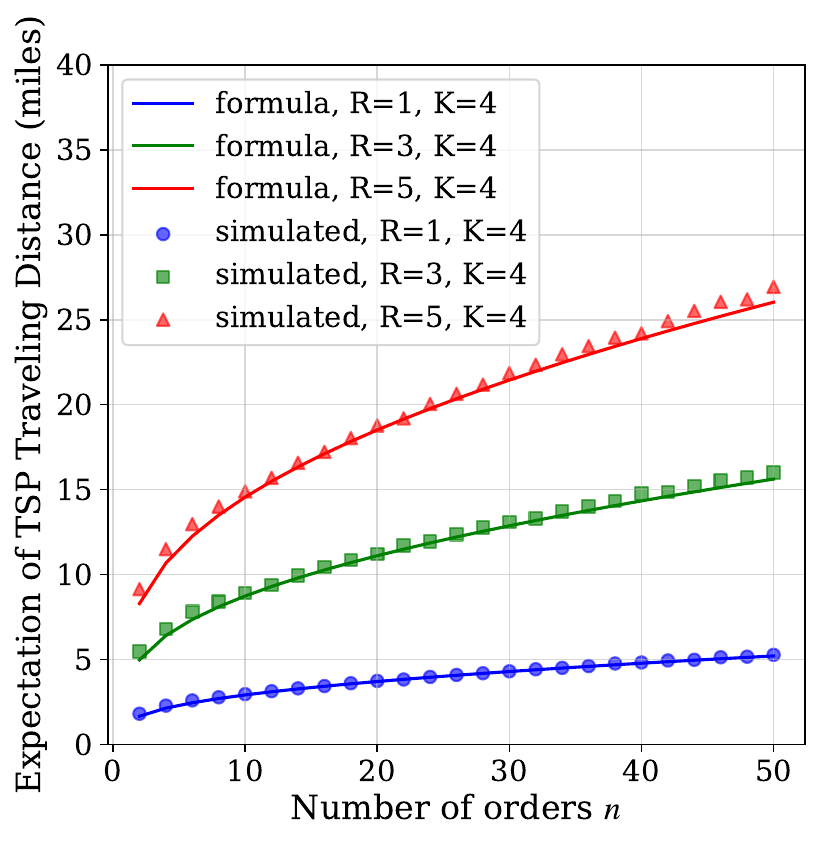}
		\caption{Mean w.r.t. $n$}
	\end{subfigure}
	\begin{subfigure}[b]{.33\textwidth}
		\centering
		\includegraphics[width=0.95\linewidth]{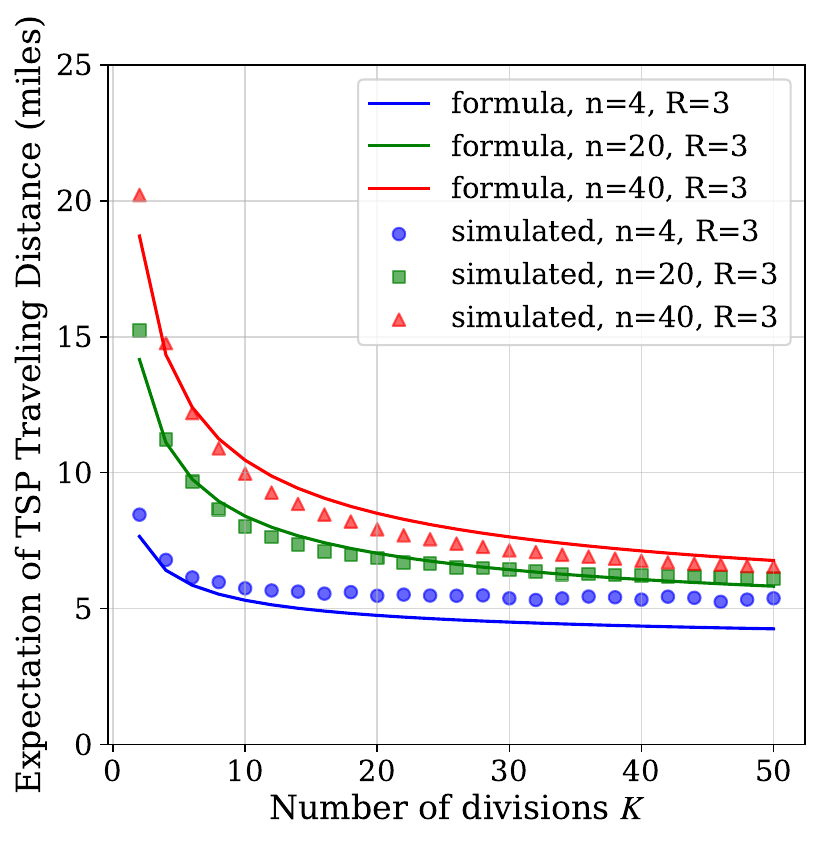}
		\caption{Mean w.r.t. $K$}
	\end{subfigure}
	% \caption{Comparison between measured and predicted expectation of the traveling distance in TSP tours w.r.t. (a) radius size $R$, (b) number of points $n$, and (c) number of divisions $K$.}
	\centering
	\begin{subfigure}[b]{.33\textwidth}
		\centering
		\includegraphics[width=0.95\linewidth]{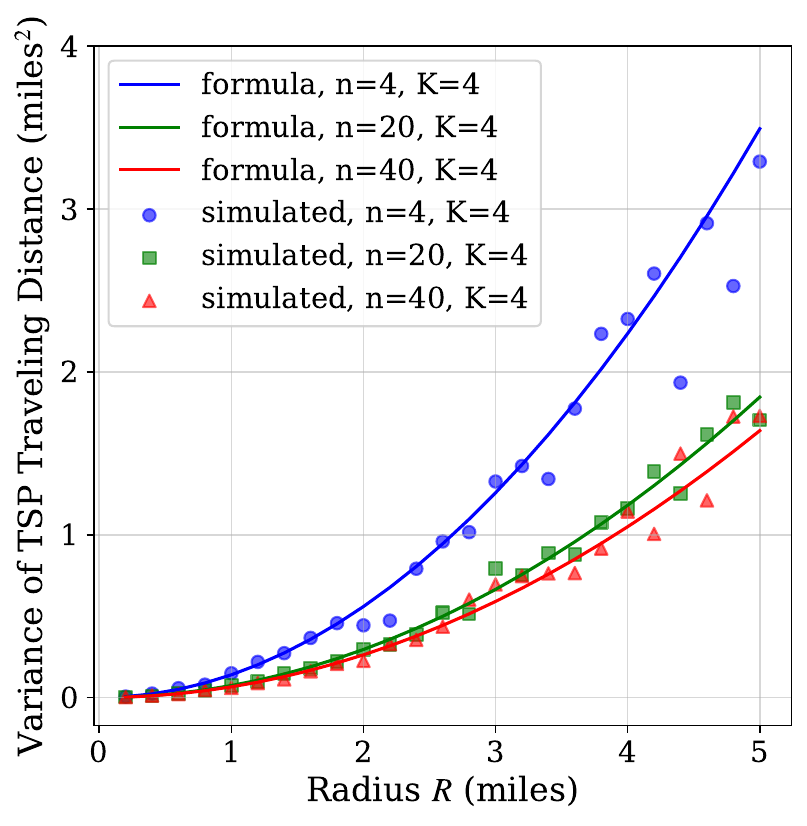}
		\caption{Var w.r.t. $R$}
	\end{subfigure}%
	\begin{subfigure}[b]{.33\textwidth}
		\centering
		\includegraphics[width=0.95\linewidth]{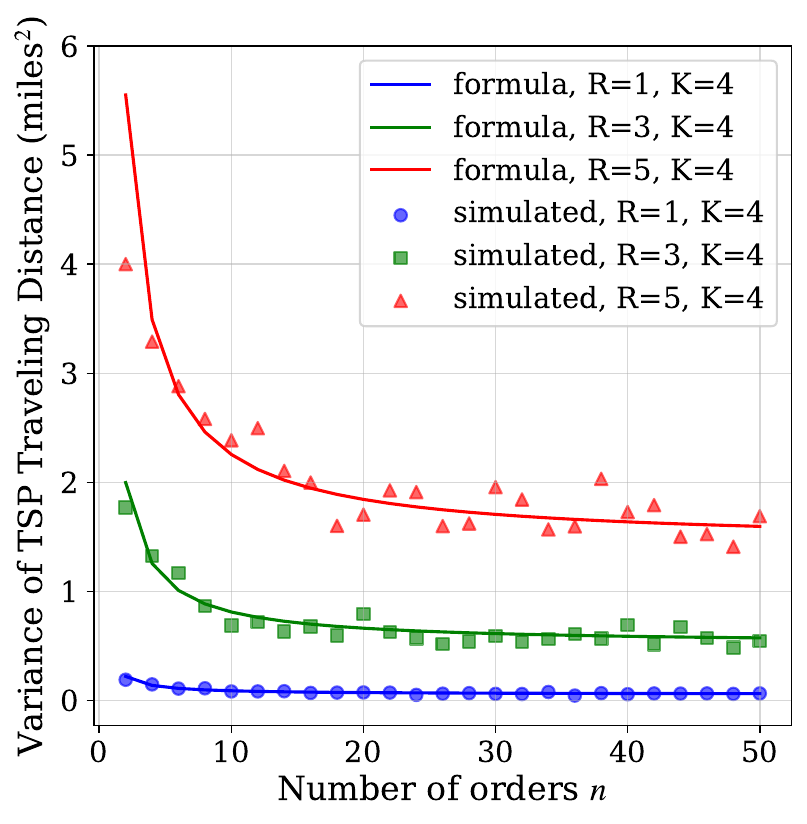}
		\caption{Var w.r.t. $n$}
	\end{subfigure}
	\begin{subfigure}[b]{.33\textwidth}
		\centering
		\includegraphics[width=0.95\linewidth]{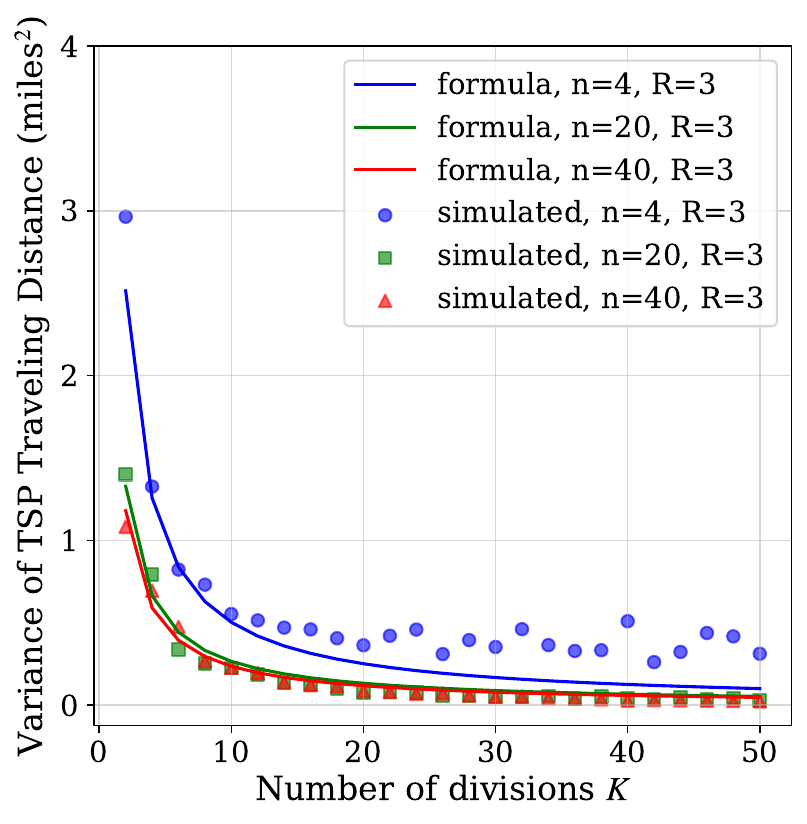}
		\caption{Var w.r.t. $K$}
	\end{subfigure}
	\caption{Comparison between measured and predicted mean and variance of the traveling distance in TSP tours w.r.t. 1) radius size $R$, 2) number of points $n$, and 3) number of divisions $K$.}
	\label{fig_TSP}
\end{figure}

We vary the values of $R$, $n$, and $K$ and simulate 1,000 trips for each case, where $n$ access points are randomly distributed across a pie-shaped circular sector with radius $R$ and angle $\theta = 2\pi / K$. The simulated results are obtained through the built-in solvers provided by Google OR-Tool's vehicle routing package, and we then calculate the average and variance of travel distance across the 1,000 simulations for each combination of $R$, $n$, and $K$. We calibrate Eq. \eqref{eq_E_Dk} and Eq. \eqref{eq_Var_Dk} for $\E[S_k]$ and $\Var[S_k]$ using the simulation results and obtain the parameter values $a = 0.64$, $b = 1.28$, $\alpha = 0.42$, and $\beta = 0.07$. Figure \ref{fig_TSP} illustrates the fitting curves for the expected value and variance of the traveling distance in TSP tours. The results indicate that the approximations using Eqs. \eqref{eq_E_Dk} and \eqref{eq_Var_Dk} achieve $R^2$ values of 98.64\% and 85.98\%, respectively, when compared to the optimal solutions provided by Google OR-Tool.

\hspace{-1.5em}  \textit{Measurement of approximations for customers' waiting time}

\begin{figure}[p!]
	\begin{subfigure}[b]{.45\textwidth}
		\centering
		\includegraphics[width=0.83\linewidth]{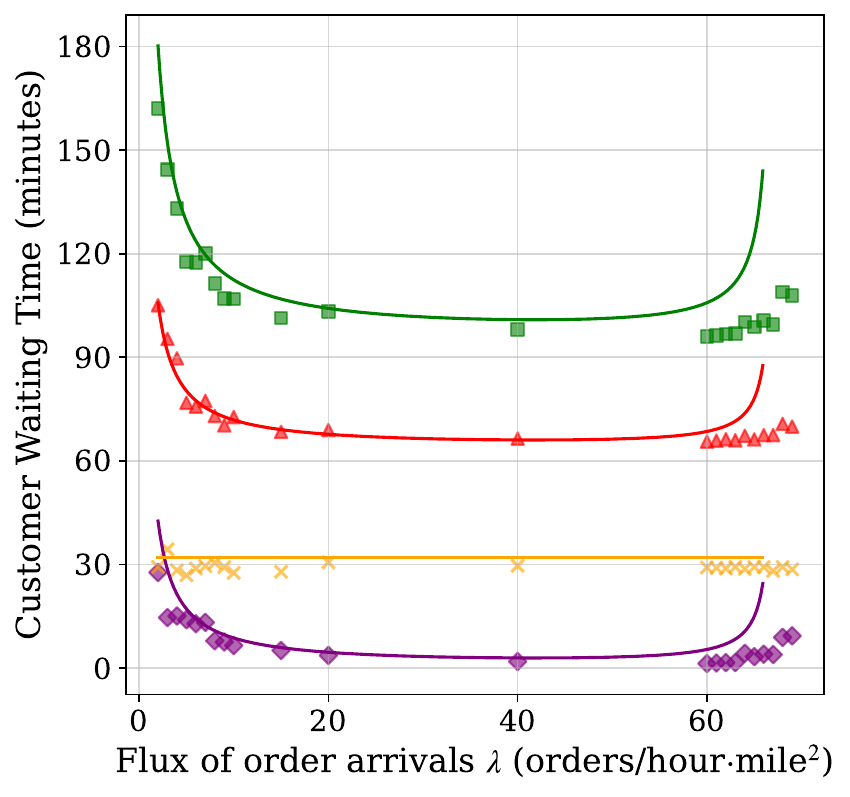}
		\caption{}
	\end{subfigure}\hfill
	\begin{subfigure}[b]{.45\textwidth}
		\centering
		\includegraphics[width=0.83\linewidth]{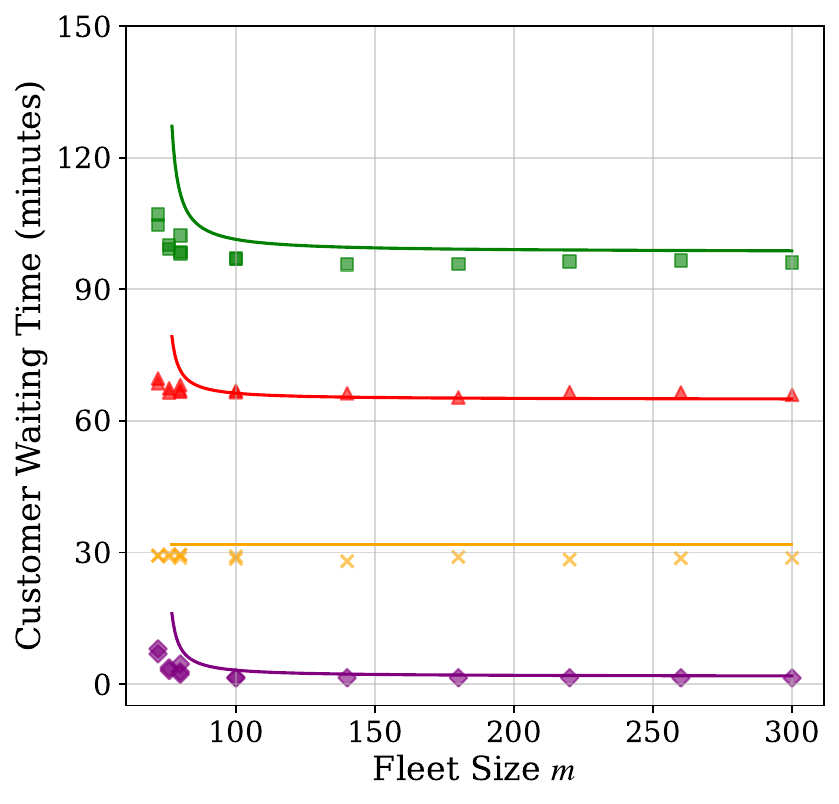}
		\caption{}
	\end{subfigure}
	\begin{subfigure}[b]{.45\textwidth}
		\centering
		\includegraphics[width=0.83\linewidth]{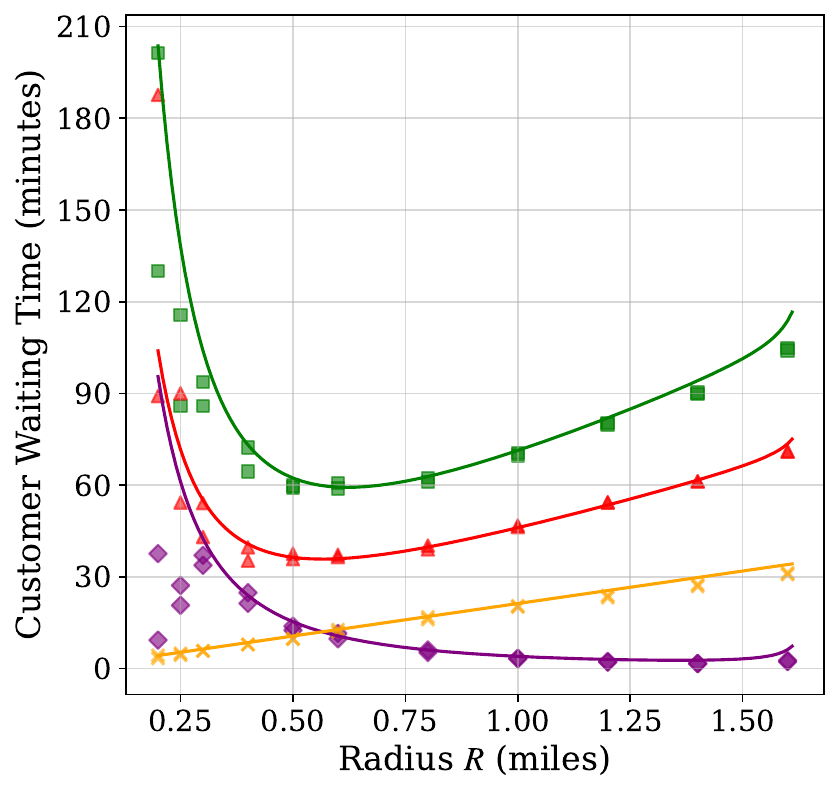}
		\caption{}
	\end{subfigure}\hfill
	\begin{subfigure}[b]{.45\textwidth}
		\centering
		\includegraphics[width=0.83\linewidth]{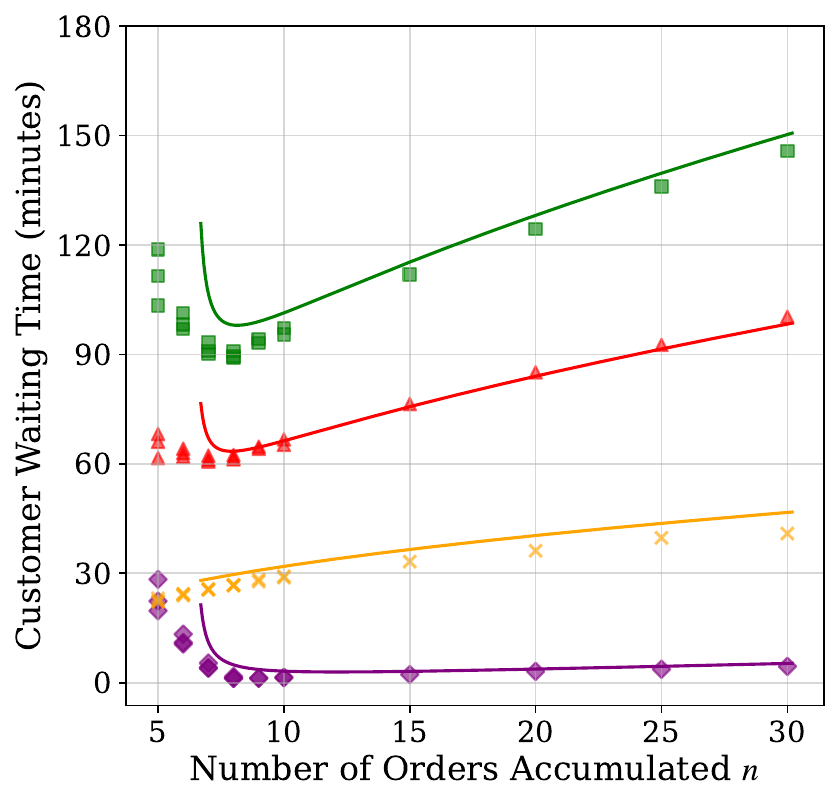}
		\caption{}
	\end{subfigure}
    \begin{center}
    	\begin{subfigure}[b]{0.9\textwidth}
    		\centering
    		\includegraphics[width=0.83\linewidth]{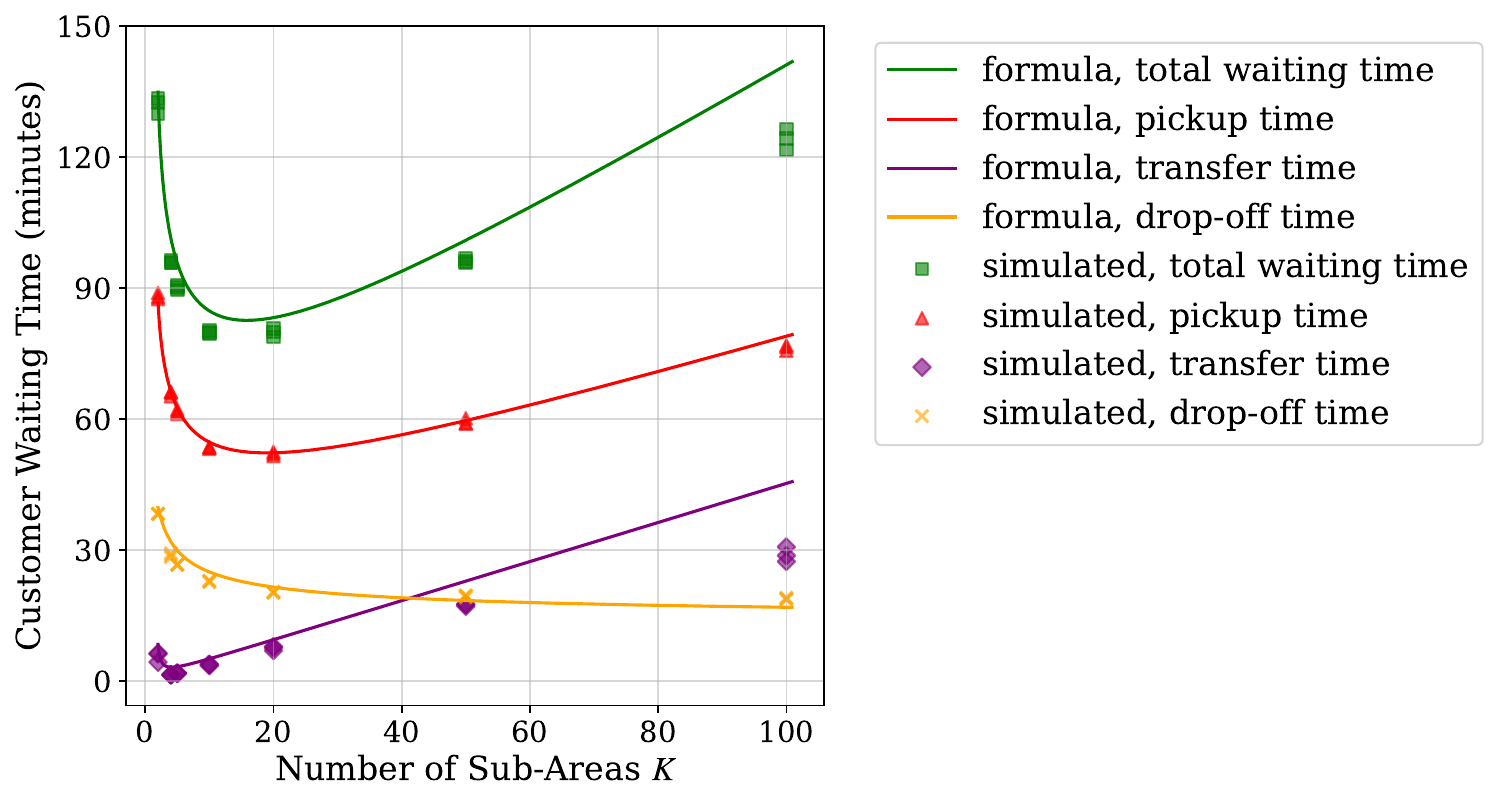}
    		\caption{}
    	\end{subfigure}\hfill
    \end{center}
	\caption{Comparison between measured and predicted customers' waiting time w.r.t. (a) order arrival flux $\lambda$, (b) fleet size $m$; (c) area radius $R$, (d) batch size of orders $n$, (e) number of sub-area partitions $K$.}
	\label{fig_measurement}
\end{figure}

The ``ground truth'' arrivals of delivery requests and customers' waiting time are simulated over a 6-hour period to experiment with the order fulfillment process. The arrival of requests follows a spatio-temporal Poisson process with a given flux $\lambda$. Once a sub-area accumulates $n$ orders and there is a deliverer available, one deliverer is dispatched to visit the $n$ access points following the routing provided by the Google OR-Tool. The simulation platform records the critical time intervals for each order, including the total waiting time from order generation to delivery, along with the three decomposed components described in Section \ref{sec:order_lifetime}. To ensure steady-state conditions, the first hour is treated as a warm-up period and excluded from the analysis. The average waiting time of all orders is calculated over the remaining five hours. The baseline parameters are set as follows: order arrival flux $\lambda=50$ orders/hour$\cdot$mi$^2$, fleet size $m=100$, area radius $R=1.5$ miles, batch size $n=10$ orders, and number of sub-area partition $K=4$. For each plot in Figure \ref{fig_measurement}, four out of the five parameters are held constant, while the total and three decomposed components of customer waiting time are plotted as functions of the varying parameter.

Figures \ref{fig_measurement}a and \ref{fig_measurement}c exhibit a non-monotonic trend in customer waiting time as a function of demand intensity and service area size, respectively. As the demand intensity $\lambda$ or service radius $R$ increases, customer waiting time initially decreases and then starts to rise. The initial decrease occurs because, when the order arrival rate is low or the service area is small, it takes longer to accumulate $n$ orders for batch service, resulting in excessively long customer wait times. Increasing the demand intensity or service area size raises the order arrival rate, thereby reducing the holding time in delivery. However, in the increasing segment, when the arrival rate becomes sufficiently high and approaches the service capacity of the available fleet, the system becomes congested, leading to longer waiting times. In addition to demand-side influences, Figure \ref{fig_measurement}b illustrates that a shortage of deliverers also results in system congestion and prolonged customer waiting times. However, once the fleet size reaches a sufficient level, adding more deliverers does not yield noticeable improvements, as many will remain idle. 

Interestingly, Figure \ref{fig_measurement}d demonstrates that a higher batch size of orders $n$ consistently results in longer customer waiting times due to the increased touring time per delivery cycle. However, in terms of VMT, as shown in Eq. \eqref{eq_miles_traveled}, batching more orders per delivery tour reduces the overall VMT by deliverers, highlighting a clear trade-off in system design. Regarding the number of partitioned sub-areas, Eqs. \eqref{eq_W_a} and \eqref{eq_W_q} indicate that finer partitioning of sub-areas (i.e., higher $K$) increases the time required to accumulate $n$ orders within each sub-area but reduces the travel distance per delivery cycle. Consequently, the total waiting time exhibits a non-monotonic trend, as depicted in Figure \ref{fig_measurement}e. These trade-offs over the choices of $n$ and $K$ justify the existence of optimal system designs with moderate configurations, rather than extreme values.

\subsection{Effectiveness of transshipment vs. non-transshipment}\label{sec:trans_vs_nontrans}

We next compare the performance of meal delivery systems with and without implementing transshipment strategies. Following the values recommended by \cite{nourbakhsh2012structured}, we set the operational cost per vehicle mile at $\pi_Q = \$2$/veh-mile and customers' value of time at $\pi_W = \$20$/hour. We analyze three scenarios where the radius of the service region is set at $R=\{1.2, 1.5, 1.8\}$ miles, and the order arrival flux $\lambda$ ranges from 5 to 170 orders/hour$\cdot$mi$^2$, respectively. For each scenario, the total number of available deliverers is set to $m=\{64, 100, 144\}$, maintaining a consistent ratio between the number of orders and deliverers. For all cases, given the area $A=\pi R^2$, fleet size $m$, and order arrival flux $\lambda$, we solve problems \eqref{opt_form_simplified} and \eqref{opt_VRPPD} to determine the optimal designs for transshipment and non-transshipment operations, respectively. We then compare the average customer waiting time (minutes/customer) and the average deliverer VMT (miles/hour-vehicle) achieved under transshipment with the corresponding metrics in the non-transshipping VRPPD counterpart.

\begin{figure}[b!]
	\centering
	\begin{subfigure}[b]{.5\textwidth}
		\centering
		\includegraphics[width=0.95\linewidth]{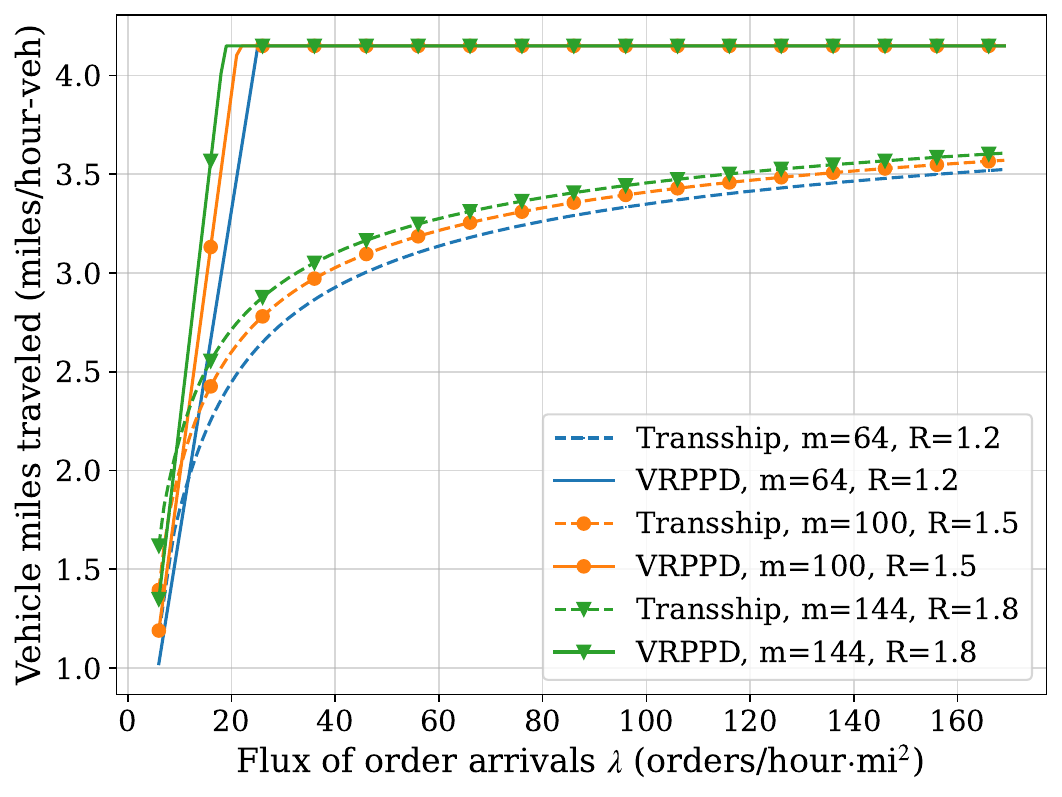}
		\caption{}
        \label{fig_optimization_lambda_a}
	\end{subfigure}%
	\begin{subfigure}[b]{.5\textwidth}
		\centering
		\includegraphics[width=0.96\linewidth]{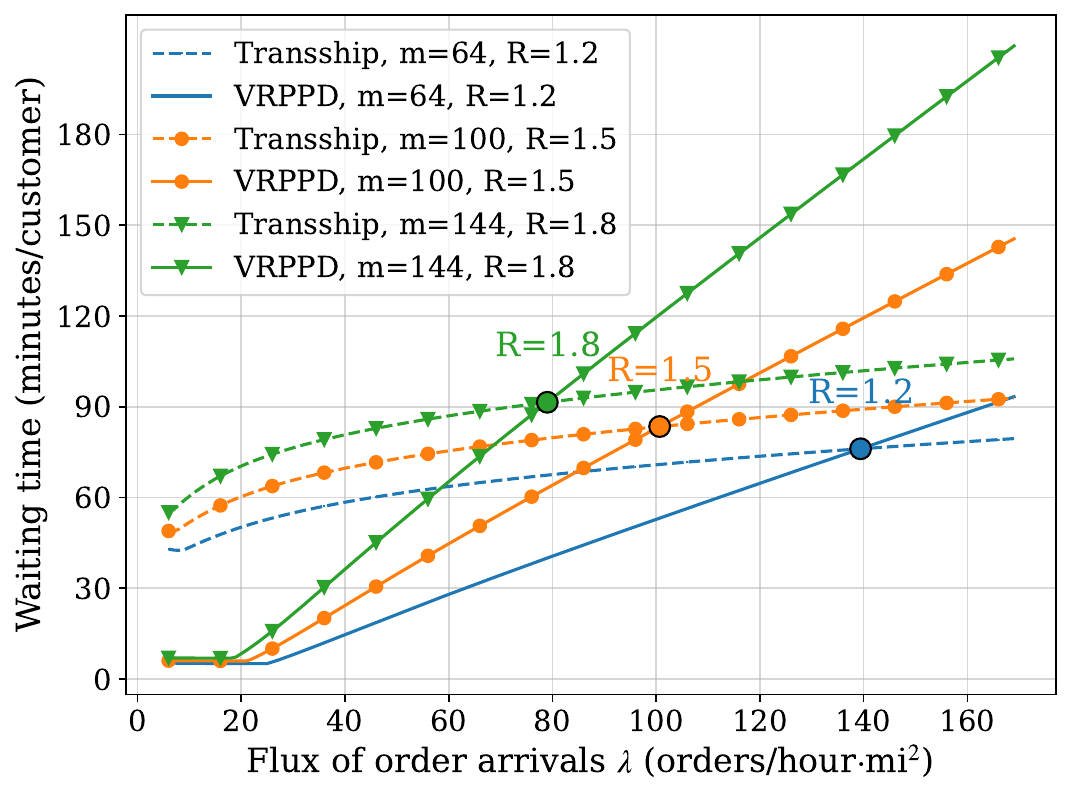}
		\caption{}
        \label{fig_optimization_lambda_b}
	\end{subfigure}
	\caption{Comparison between the transshipment and non-transshipment strategies on (a) VMT and (b) customer waiting time w.r.t. order arrival flux.}
	\label{fig_optimization_lambda}
\end{figure}

\begin{figure}[th!]
	\centering
	\includegraphics[width=0.9\textwidth]{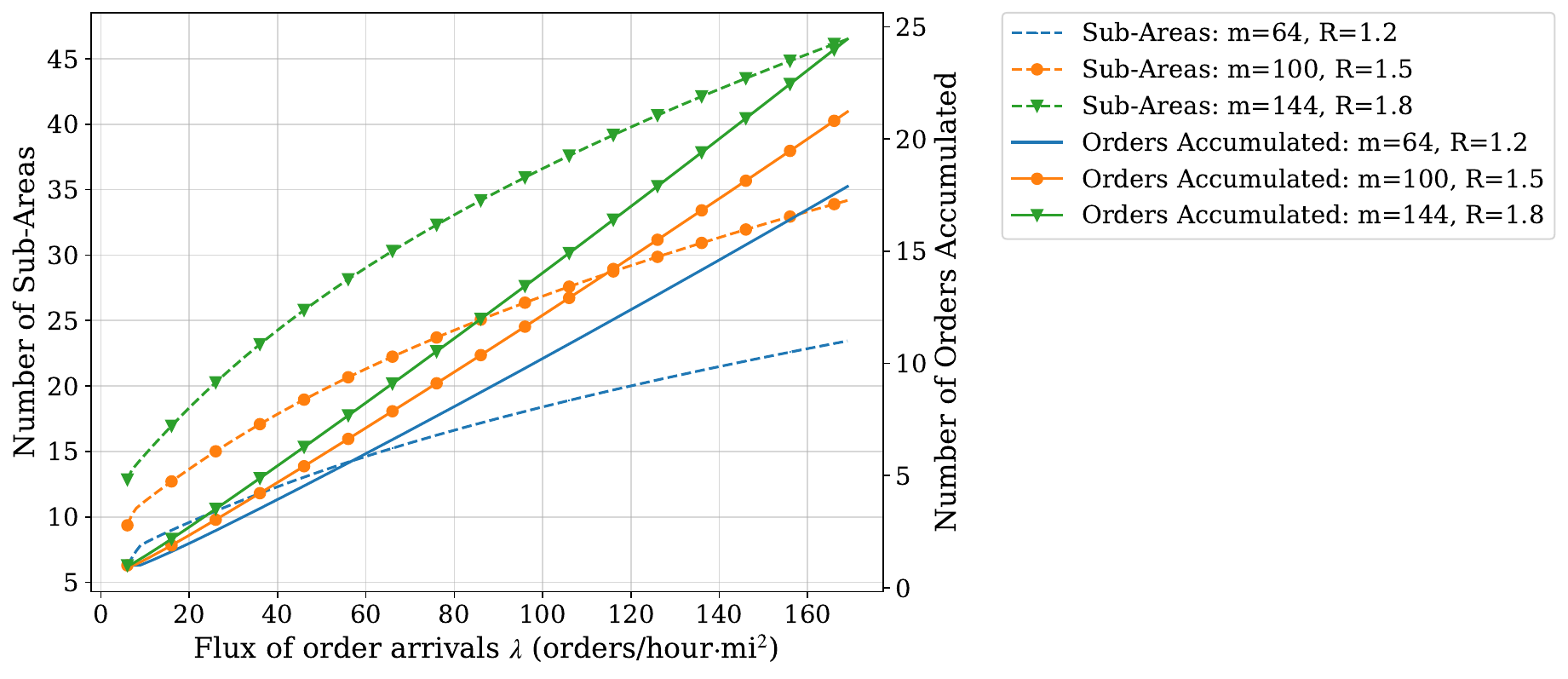}
	\caption{Optimal partitioning and batching strategic decisions for meal delivery system with transshipment w.r.t. order arrival flux.}
	\label{fig_optimization_K_n}
\end{figure}

Figure \ref{fig_optimization_lambda} presents a comparison of results and highlights the conditions under which transshipment demonstrates relative strength over non-transshipment strategies. From Figure \ref{fig_optimization_lambda_a}, it is evident that while both strategies involve order consolidation onto individual vehicles, the localized batching facilitated by transshipment consistently reduces VMT compared to the VRPPD approach. Without transshipment, when demand exceeds a certain threshold---around $\lambda \approx 20$ in Figure \ref{fig_optimization_lambda_a}---all deliverers are fully deployed and travel continuously between pickup and drop-off points, causing the average VMT per hour to converge toward the vehicle's operating speed. In contrast, our transshipment design allows deliverers to wait at the microhub until the number of visiting points accumulates to the batching size $n$, thereby resulting in significant VMT savings. 

Furthermore, as shown in Figure \ref{fig_optimization_lambda_b}, for a given service range $R$ with a fixed supply, the transshipment strategy demonstrates its advantages by achieving lower customer waiting times compared to the VRPPD counterparts as demand increases within the examined region. In Figure \ref{fig_optimization_K_n}, we demonstrate the optimal partitioning and batching results obtained from solving the optimal design problem \eqref{opt_form_simplified}. It suggests that, with sufficient service capacity, system efficiency can potentially improve as demand grows by raising the number of sub-areas and the batching size of orders. However, at lower demand levels, transshipment through the microhub becomes less efficient than the direct pickup-and-delivery strategy due to the additional time required for routing and transferring. Notably, Figure \ref{fig_optimization_lambda_b} also illustrates that the efficiency gains enabled by transshipment are influenced by the size of the service region. While keeping the ratio of total demand to total supply constant, larger service areas tend to exacerbate the challenges faced by the VRPPD model due to the increased frequencies of long-distance cross-regional delivery tours. In contrast, the localization strategy enabled by transshipment mitigates the impact of regional size on meal delivery and maintains efficiency by dynamically adjusting the number of sub-areas to accommodate varying sizes of the service region. In summary, transshipment through a microhub offers significant advantages over non-transshipment strategies in scenarios with high demand or insufficient supply, particularly in servicing areas of larger sizes.

\subsection{A case study of Meituan meal delivery}\label{sec:case_study}

\begin{figure}[t!]
	\centering
	\begin{subfigure}[b]{.5\textwidth}
		\centering
		\includegraphics[width=0.95\linewidth]{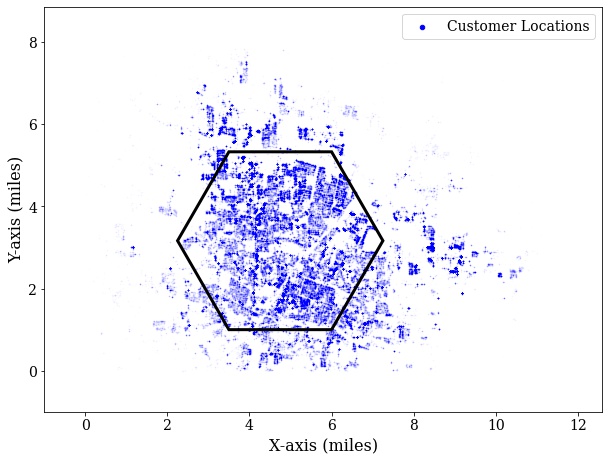}
		\caption{}
        \label{fig_Meituan_1a}
	\end{subfigure}%
	\begin{subfigure}[b]{.5\textwidth}
		\centering
		\includegraphics[width=0.95\linewidth]{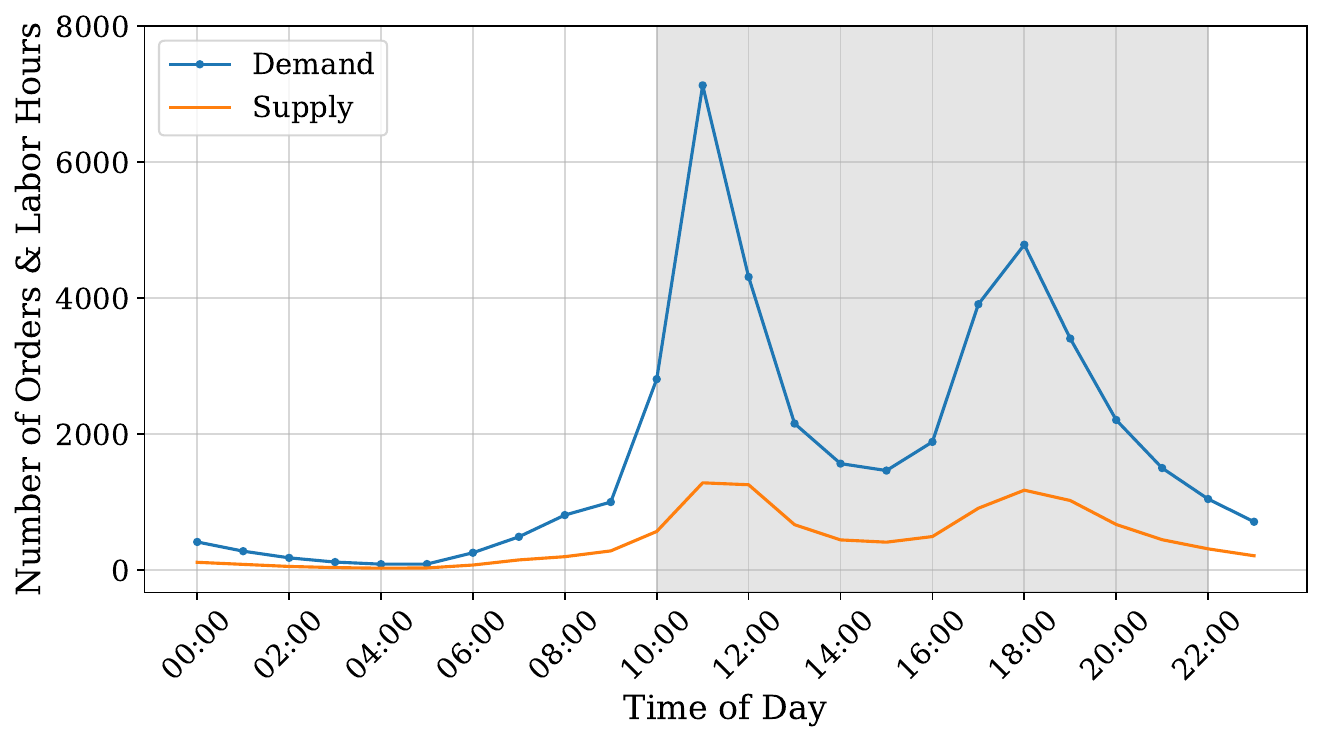}
		\caption{}
        \label{fig_Meituan_1b}
	\end{subfigure}
	\caption{Visualization of Meituan meal delivery data: (a) a map showing customer locations and the targeted area with a radius of $R=2.5$ miles for analysis; (b) hourly variations in demand and supply within the targeted service area throughout the day.}
	\label{fig_Meituan_1}
\end{figure}

We conduct a case study comparing transshipment and non-transshipment models using Meituan meal delivery data. The original dataset comprises 654,343 order requests recorded between October 17 and 24, 2022, in an undisclosed city, publicly provided by Meituan. After filtering out requests located far from the central cluster, the master dataset retains 484,650 order records, 4,571 deliverers, and 206,748 trajectories. In Figure \ref{fig_Meituan_1a}, we define a hexagonal service area with a radius of $R=2.5$ miles, intended for transshipment operations through a single microhub, and assume that orders are uniformly and randomly distributed within the hexagon. For Figure \ref{fig_Meituan_1b}, demand is defined as the number of orders occurring per hour within the service area, while supply is measured as the total active labor hours of deliverers derived from their trajectory data. Both demand and supply values are averaged over the week to represent a typical day. 

\begin{figure}[t!]
	\centering
	\includegraphics[width=1\textwidth]{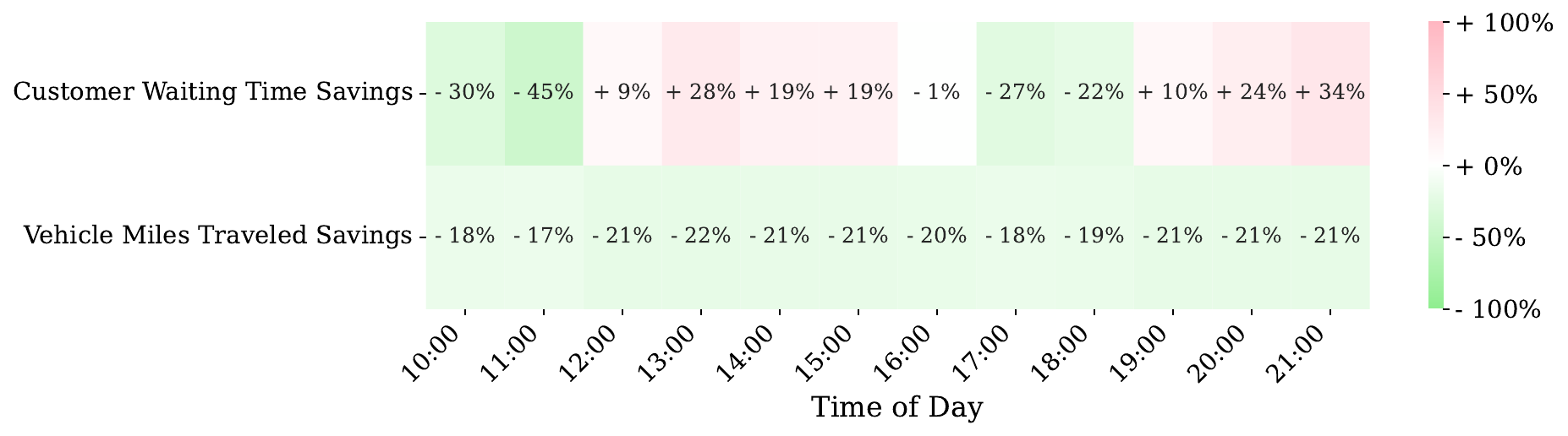}
	\caption{Savings in deliverers' vehicle miles traveled and customers' average waiting time achieved by transshipment compared to non-transshipment.}
	\label{fig_Meituan_demand_supply_1}
\end{figure}

For each pair of demand and supply data from 10:00 to 21:00, we independently solve the optimal design problem with transshipment \eqref{opt_form_simplified} and without \eqref{opt_VRPPD}. The relative differences in two key metrics---customers' average waiting time and deliverers' vehicle miles traveled---are shown in Figure \ref{fig_Meituan_demand_supply_1}, where greener shades indicate that the transshipment model is more efficient than the non-transshipment counterpart. As anticipated from previous analysis, the transshipment strategy constantly reduces VMT across all time periods compared to the non-transshipment approach. For customer waiting time, during lunch and dinner peak hours, when demand significantly exceeds supply, the transshipment strategy achieves up to a 45\% reduction, with the greatest saving occurring between 11:00 and 12:00, coinciding with the highest demand period of the day. However, during non-peak hours, the VRPPD model without transshipment outperforms the transshipment strategy in providing shorter customer waiting times. These findings suggest an alternating approach for meal delivery, where transshipment is adopted during peak hours, while the system reverts to the classic pickup-and-delivery model during non-peak hours. 

In Appendix \ref{sec:appd_Meituan}, we reduce the service range covered by individual microhubs for transshipment and deploy multiple microhubs to cover the entire region. The results further confirm that transshipment is particularly advantageous in high-demand areas during peak hours.

\section{Conclusion}\label{sec:conclusion}

This paper investigates a novel meal delivery strategy that leverages a microhub for transshipment. An analytical model was developed using a CA method to evaluate the system's expected performance under the transshipment strategy. To estimate customer waiting time, we decomposed the lifecycle of a meal order and adapted a queuing network to approximate the time spent by orders at different service stages. A comprehensive set of simulations was conducted to validate the accuracy of these approximations. The analytical approximations were then incorporated into an optimal design problem to solve for the performance envelope under various exogenous market conditions. The numerical results indicate that when the system is optimally configured, the proposed meal delivery strategy through transshipment significantly reduces vehicle miles traveled compared to the conventional VRPPD strategy without transshipment. Additionally, the transshipment strategy effectively decreases expected customer waiting time, particularly under high demand or low supply conditions, and proves especially advantageous for servicing larger areas. However, in scenarios with sufficient supply, excessive reliance on transshipment may lead to longer customer waiting times. These insights into the applicability of the transshipment strategy were finally validated using Meituan meal delivery data.

Future research can explore several promising directions. First, this study examines the simplest scenario of a single, centrally placed microhub for transshipment in an isotropic area. A natural extension would be to analyze heterogeneous settings where both demand and supply are unevenly distributed and determine the optimal location of the microhub to accommodate spatial variations in delivery patterns. Second, future research can also explore scenarios involving multiple microhubs optimally positioned to facilitate transshipment within a larger area. In such cases, a two-echelon delivery system could be investigated, where micromobility options such as e-bikes or e-scooters are dedicated to servicing smaller areas covered by individual microhubs, while heavy-duty delivery vehicles periodically transport packages between microhubs. Third, toward the end of the paper, we sketchily rationalize an alternating strategy for meal delivery to adopt transshipment during peak periods and switch back to non-transshipment strategies during non-peak periods. However, the current modeling framework treats different time periods independently, without factoring in transition costs, start-up delays, or relocation complexities. Future work could incorporate these dynamic transition effects to provide a more realistic evaluation of the dynamic system design. Lastly, it would be valuable to endogenize both customer demand and deliverer supply in the model. Customer demand might fluctuate in response to expected waiting time, while deliverer availability could be influenced by overall travel distances and order completion rates. Integrating these behavioral and economic factors would yield more accurate predictions of system performance and help platform operators design more robust and adaptive delivery strategies.

\section*{Acknowledgement}
This research utilized data from the TSL Data-Driven Research Challenge, provided by Meituan. We extend our heartfelt gratitude to Meituan for sharing their service data. This work is partially supported by the National Science Foundation (CMMI-2443338).

\section*{Declaration of generative Al and Al-assisted technologies in the writing process}
During the preparation of this work the authors used ChatGPT in order to enhance the readibility of the paper. After using this tool/service, the authors reviewed and edited the content as needed and took full responsibility for the content of the published article.

\newpage
\section*{Nomenclature}
\LTcapwidth=\textwidth
\begin{table}[ht!]
	\centering
	\caption{Notation list of parameters and variables}
	\label{notatoin1}
	\begin{tabular}{c l c} 
	\hline
	Notation & Description & Unit\\
	\hline\hline
	\multicolumn{3}{l}{\textbf{\emph{Input Parameters}}}\\
        $\lambda$ & Expected number of orders generated per hour per unit area & 1/hour$\cdot$mi$^2$ \\
        $m$ & Total size of delivery fleet & \\
        $v$  & Speed of delivery vehicles & mph\\
        $A$  & Total area of the region & mi$^2$ \\
        $a, b, \alpha, \beta$ & Parameters related to TSP approximations in a pie-shaped area &  \\
        $\pi_Q$  & Operating cost per vehicle mile traveled & \$ \\
        $\pi_W$  & Value of time & \$/hour\\
	\hline
	
	\multicolumn{3}{l}{\textbf{\emph{Decision Variables}}}\\
        $K$ & Number of sub-area partitions & \\
        $A_k$  & Area of each sub-area $k$ & mi$^2$\\
        $\lambda_k$ & Expected number of new orders generated per hour in sub-area $k$ & 1/hour\\
        $\delta_k$  & Number of pickup and drop-off locations to visit in sub-area $k$ & 1/hour\\
        $p_{i,k}$ & Probability of an order picked up in sub-area $i$ being dropped off in $k$ \\
        $m_k$ & Number of deliverers assigned to sub-area $k$ \\
        $n_k$ & Batch size of packages per delivery tour in sub-area $k$ \\
        % $D_k$  & Traveling distance of individual deliverer in sub-area $k$ per delivery cycle & mi\\
        $S_k$  & Traveling time of individual deliverer in sub-area $k$ per delivery cycle & hour\\
        $Q$  & Expected vehicle miles traveled by deliverers per hour of operation & mph \\
        $W_k^a$  & Average waiting time for accumulating $n_k$ packages in sub-area $k$ & hour\\
        $W_k^q$  & Average package holding time for available deliverers in sub-area $k$ & hour\\
        $W$  & Average total waiting time per customer & hour\\
	\hline
	\end{tabular}
\end{table}
\clearpage

\begin{appendices}

\section{Proof of Propositions} \label{sec:appd_proof}
\setcounter{equation}{0}
\setcounter{figure}{0}
\setcounter{table}{0}
\renewcommand\theequation{A\arabic{equation}}
\renewcommand\thefigure{A\arabic{figure}}
\renewcommand\thetable{A\arabic{table}}

This appendix presents the proofs for the propositions introduced in \Cref{sec:microhub,sec:VRPPD}.

\subsection{Proof of Proposition \ref{proposition_1}}

Consider a pie-shaped circular sector with radius $R$ and central angle $\theta$, containing $N$ nodes randomly distributed across the entire region. Let $D_i$ denote the distance from point $i$ within the sector to the center of the circle. Then, the distance from the origin to the farthest node is defined as $R' = \max\{D_i\}_{i\in N}$, and its CDF and PDF are given by the following expressions:
\begin{align*}
    F_{R'}(r) &= \probP(\max\{D_i\}_{i\in N} \leq r) = \probP(D_1 \leq r, \cdots, D_N \leq r) = \left( \frac{\frac{\theta}{2\pi} \pi r^2}{\frac{\theta}{2\pi} \pi R^2} \right)^N = \frac{r^{2N}}{R^{2N}}, \\
    f_{R'}(r) &= 2N\frac{r^{2N-1}}{R^{2N}}.
\end{align*}
The expected value of $R'$ is thus given by
\begin{align*}
    \E[R'] &= \int_0^R r \cdot f_{R'}(r)\ \d r = \frac{2N}{2N+1} R.
\end{align*}

\subsection{Proof of Proposition \ref{proposition_2}}

For the Spatial Poisson process, let $\delta=\frac{N_O+N_D}{\pi R^2}$ denote the density of access points across the entire area. The survival function and PDF of the distance $X$ from a point to its nearest neighboring point are given by
\begin{align*}
    & \probP(X>x) = \probP(\text{no points within distance $x$}) = e^{-\delta \cdot \pi x^2}, \\
    & f_X(x) = \delta \cdot 2\pi x \cdot e^{-\delta \cdot \pi x^2}, \ x \geq 0.
\end{align*}
Similarly, for the variable  $Y$, which follows the Rayleigh distribution, the survival function and PDF are given by:
\begin{align*}
    & \probP(Y>y) = e^{-\frac{y^2}{2\sigma^2}}, \\
    & f_Y(y) = \frac{y}{\sigma^2}\cdot e^{-\frac{y^2}{2\sigma^2}}, \ y \geq 0.
\end{align*}
Thus, the probability of $X>Y$ can be derived as:
\begin{align*}
    \probP(X>Y) &= \int_0^\infty \int_y^\infty f_X(x)\ f_Y(y) \, \d x\ \d y \\
    &= \int_0^\infty \int_y^\infty 2\delta \pi x e^{-\delta \pi x^2} \cdot \frac{y}{\sigma^2}e^{-\frac{y^2}{2\sigma^2}} \, \d x\ \d y\\
    &= \int_0^\infty \left( \frac{y}{\sigma^2}e^{-\frac{y^2}{2\sigma^2}} \right) \left( \int_y^\infty 2\delta \pi x e^{-\delta \pi x^2} \, \d x \right) \, \d y.
\end{align*}
First, considering the integral over $x$, we have $\probP(X>y) = \int_y^\infty 2\delta \pi x e^{-\delta \pi x^2} \, \d x =  e^{-\delta \pi y^2}$. Substituting it back yields $\probP(X>Y) = \int_0^\infty \frac{y}{\sigma^2} \cdot e^{-(\delta \pi + \frac{1}{2\sigma^2}) y^2} \, \d y$. By letting $a = \delta \pi + \frac{1}{2\sigma^2}$, $u=a y^2$, and $\d u=2ay\d y$, we obtain
\begin{align*}
    \probP(X>Y) &= \frac{1}{\sigma^2} \cdot \frac{1}{2a} \int_0^\infty e^{-u} \d u \\
    &= \frac{1}{1+2\delta \pi \sigma^2} \\
    &= \frac{1}{1+\frac{2(N_O+N_D)\sigma^2}{R^2}}.
\end{align*}

\subsection{Proof of Proposition \ref{proposition_3}}
For variable $X$, the survival function is given by:
\begin{align*}
    \probP(X>z) &= \probP(\text{no points within distance $z$}) = e^{-\delta \cdot \pi z^2} = e^{-\frac{N_O+N_D}{R^2} z^2}.
\end{align*}
For variable $Y$, by substituting $\sigma = \frac{1}{\sqrt{2\mu}}$, its survival function becomes:
\begin{align*}
    \probP(Y>z) &= e^{-\frac{z^2}{2\sigma^2}} = e^{-\mu z^2}.
\end{align*}
The expectation of $\min(X,Y)$ is then given by:
\begin{align*}
    \E[\min(X,Y)] &= \int_0^{+\infty} \probP(X>z, Y>z) \, \d z \\
    &= \int_0^{+\infty} \probP(X>z) \cdot \probP(Y>z) \, \d z \\
    % &= \int_0^{+\infty} e^{-\frac{N_O+N_D}{R^2} z^2} \cdot e^{-\mu z^2} \, \d z \\
    &= \int_0^{+\infty} e^{-\left( \frac{N_O+N_D}{R^2}+\mu \right) z^2} \, \d z \\
    & \overset{\text{Gaussian integral}}{=\joinrel=} \frac{1}{2} \sqrt{\frac{\pi}{\frac{N_O+N_D}{R^2}+\frac{1}{2\sigma^2}}}.
\end{align*}

\section{Validation of Rayleigh Distribution Assumption for Order Distances}\label{sec:appd_Rayleigh}
\setcounter{equation}{0}
\setcounter{figure}{0}
\setcounter{table}{0}
\renewcommand\theequation{B\arabic{equation}}
\renewcommand\thefigure{B\arabic{figure}}
\renewcommand\thetable{B\arabic{table}}

\begin{figure}[h!]
	\centering
	\includegraphics[width=0.7\textwidth]{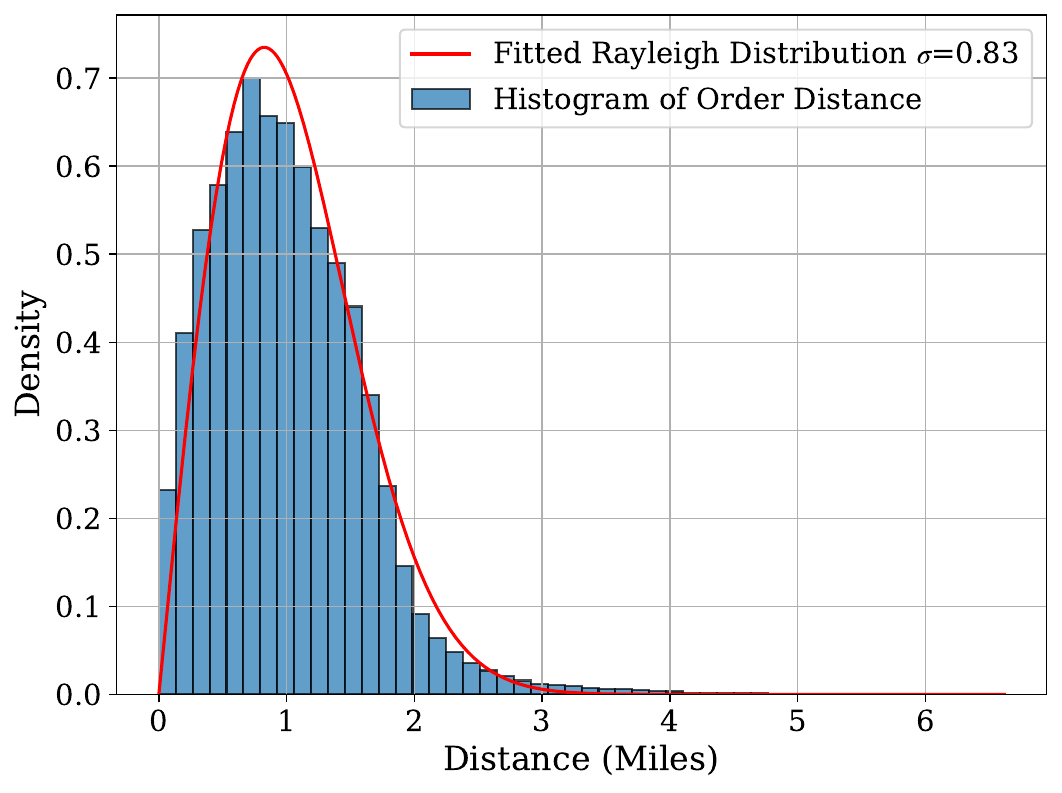}
	\caption{Fitting of the Rayleigh distribution to meal delivery order distances.}
	\label{fig_Rayleigh}
\end{figure}

This appendix details the calibration and validation of using Rayleigh distribution to model the delivery distance $d$ between the origin and destination of meal orders. The PDF of the Rayleigh distribution is given by:
\begin{align*}
    f(d| \sigma) &= \frac{d}{\sigma^2} e^{-\frac{d^2}{2\sigma^2}}, \quad d > 0, \; \sigma > 0,
\end{align*}
where $\sigma$ is the scale parameter. Given the observed distances $\{d_1, d_2, \cdots, d_n\}$ for $n$ orders, the log-likelihood function is:
\begin{align*}
    \mathcal{L}(\sigma) &= \sum_{i=1}^n \left( \ln d_i - 2\ln \sigma - \frac{d_i^2}{2\sigma^2} \right).
\end{align*}
To obtain the maximum likelihood estimate (MLE) of $\sigma$, we differentiate $\mathcal{L}(\sigma)$ with respect to $\sigma$, set the derivative to zero, and solve:
\begin{align*}
    \frac{\partial \mathcal{L}(\sigma)}{\partial \sigma} &= -\frac{2n}{\sigma} + \frac{\sum_{i=1}^n d_i^2}{\sigma^3} = 0.
\end{align*}
From the Meituan meal delivery order data, the MLE of $\sigma$ is calculated as:
\begin{align*}
    \hat{\sigma} &= \sqrt{\frac{\sum_{i=1}^n d_i^2}{2n}} = 0.83.
\end{align*}
Figure \ref{fig_Rayleigh} compares the fitted Rayleigh distribution with the histogram of 423,449 observed delivery distances for meal orders recorded during the week of October 17–24, 2022, between 10:00 AM and 10:00 PM each day. The average delivery distance is approximately 1 mile, with a higher concentration of orders occurring below this distance. This skewed distribution is well captured by the Rayleigh model, which achieves a mean squared error of 0.001249 miles.

\section{Case study with smaller transshipment range}\label{sec:appd_Meituan}
\setcounter{equation}{0}
\setcounter{figure}{0}
\setcounter{table}{0}
\renewcommand\theequation{B\arabic{equation}}
\renewcommand\thefigure{B\arabic{figure}}
\renewcommand\thetable{B\arabic{table}}

\begin{figure}[th!]
	\centering
	\includegraphics[width=0.7\textwidth]{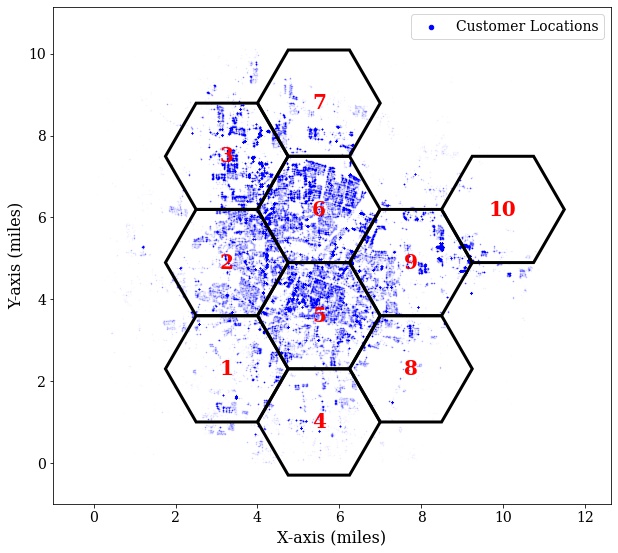}
	\caption{A map showing customer locations and the service range covered by each microhub, with a radius of $R=1.5$ miles.}
	\label{fig_Meituan_location_2}
\end{figure}

\begin{figure}[th!]
	\centering
	\includegraphics[width=1\textwidth]{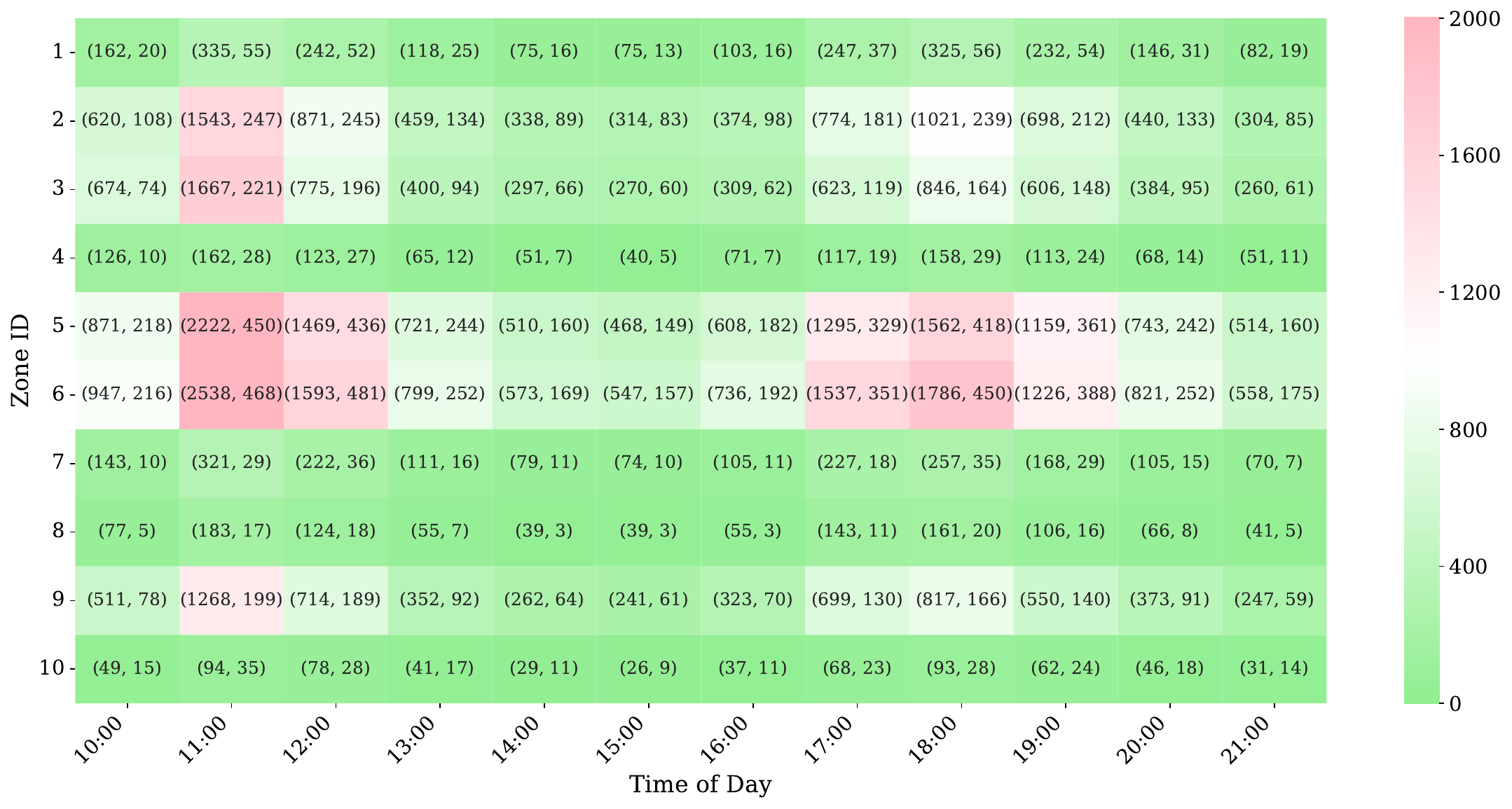}
	\caption{Hourly variations in demand and supply across the ten hexagons throughout the day.}
	\label{fig_Meituan_demand_supply_2}
\end{figure}

In this appendix, we reduce the transshipment service range through a single microhub from $R=2.5$ in \Cref{sec:case_study} to $R=1.5$, dividing the entire region is divided into 10 hexagons, each with distinct demand and supply patterns. Each microhub is used exclusively for transshipment within its respective hexagon. Figure \ref{fig_Meituan_location_2} illustrates the zonal divisions, while Figure \ref{fig_Meituan_demand_supply_2} presents the demand and supply for each hexagon at different times of the day. The two numbers in each cell indicate the number of orders and the labor hours available, respectively. Those hexagons with higher demand are highlighted in redder shades.

\begin{figure}[th!]
	\centering
	\includegraphics[width=1\textwidth]{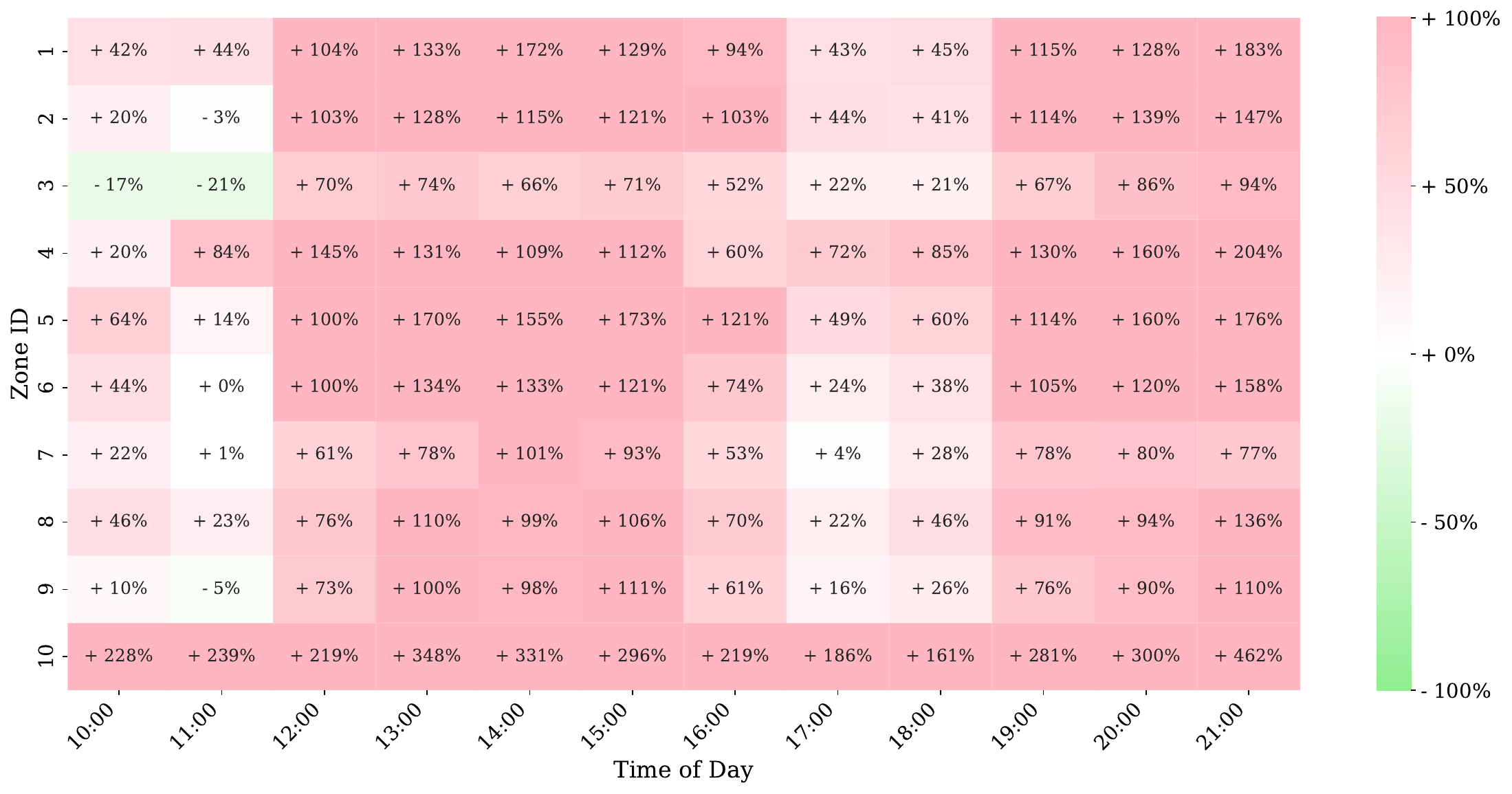}
	\caption{Savings in customer waiting time achieved through transshipment across the ten hexagonal areas throughout the day.}
	\label{fig_Meituan_demand_supply_W_2}
\end{figure}

\begin{figure}[th!]
	\centering
	\includegraphics[width=1\textwidth]{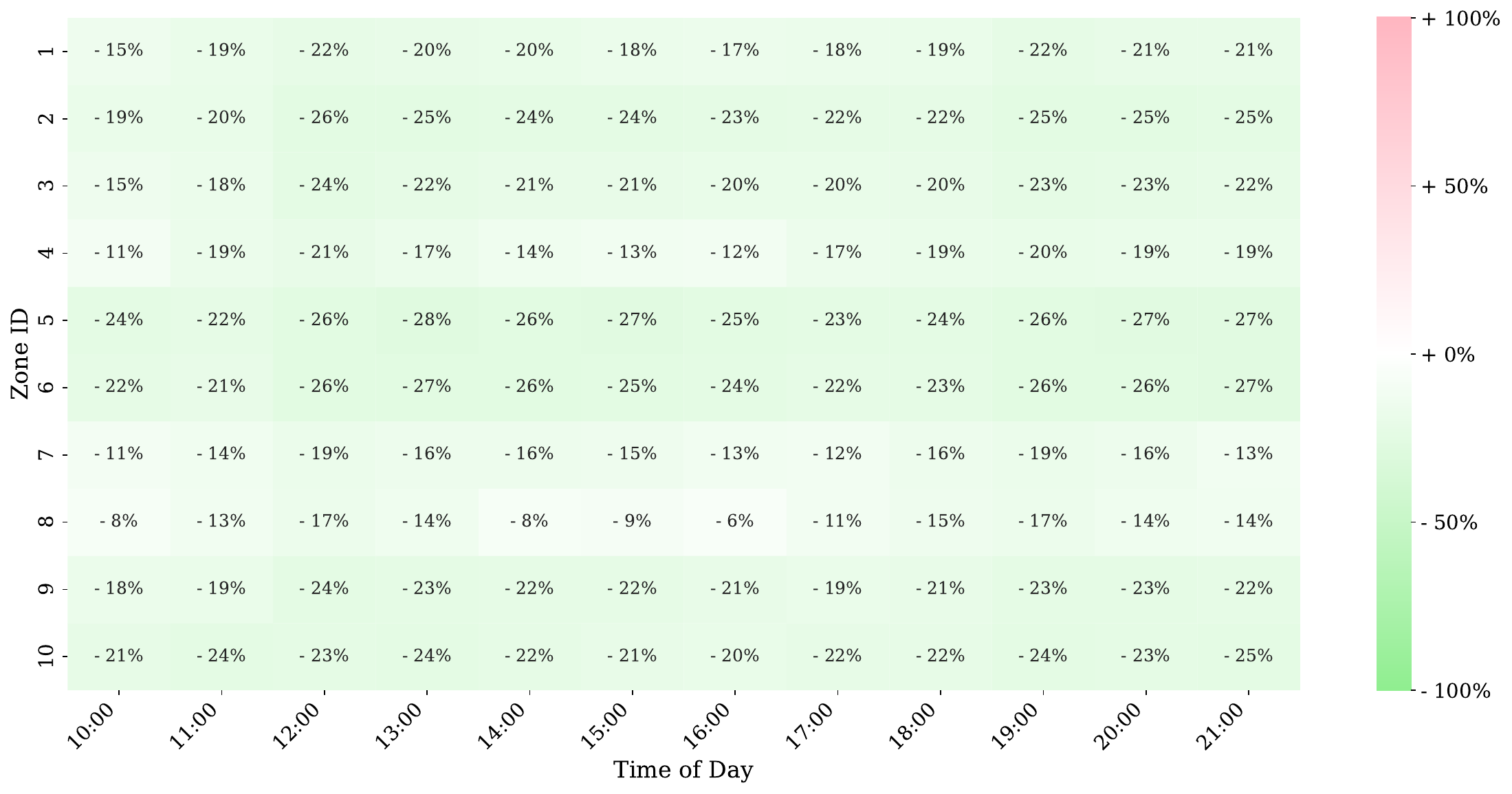}
	\caption{Savings in vehicle miles traveled by deliverers through transshipment across the ten hexagonal areas throughout the day.}
	\label{fig_Meituan_demand_supply_Q_2}
\end{figure}

Figures \ref{fig_Meituan_demand_supply_W_2} and \ref{fig_Meituan_demand_supply_Q_2} illustrate the relative savings achieved by comparing transshipment against non-transshipment in terms of average customer waiting time and vehicle miles traveled by deliverers across different areas and time periods. Among the ten divisions, Hexagons 2, 3, 5, and 6 are identified in Figure \ref{fig_Meituan_demand_supply_2} as the busiest areas. In Figure \ref{fig_Meituan_demand_supply_W_2}, Hexagons 2 and 3 demonstrate improved customer waiting times when transitioning to transshipment during the lunch peak hour, which corresponds to the highest demand-to-supply ratio. However, in Hexagons 5 and 6, despite having the highest order volumes among all zones, the available supply is also sufficient, making the transshipment strategy less effective compared to the non-transshipping counterpart. Comparing the results with those in \Cref{sec:case_study} suggests that transshipment becomes less efficient on smaller service ranges, which aligns with the findings in \Cref{sec:trans_vs_nontrans}. While VMT savings are consistently observed in all scenarios when using transshipment, the benefits in customer waiting time remain highly sensitive to the balance between demand and supply.

\clearpage
\end{appendices}

\bibliographystyle{apalike} %dcu
%\nocite{*}
\bibliography{mybib}

\begin{thebibliography}{}

\bibitem[Agatz et~al., 2024]{agatz2024transportation}
Agatz, N., Cho, S.-H., Sun, H., and Wang, H. (2024).
\newblock Transportation-enabled services: Concept, framework, and research opportunities.
\newblock {\em Service Science}, 16(1):1--21.

\bibitem[Ansari et~al., 2018]{ansari2018advancements}
Ansari, S., Ba{\c{s}}dere, M., Li, X., Ouyang, Y., and Smilowitz, K. (2018).
\newblock Advancements in continuous approximation models for logistics and transportation systems: 1996--2016.
\newblock {\em Transportation Research Part B: Methodological}, 107:229--252.

\bibitem[Campbell, 2013]{campbell2013continuous}
Campbell, J.~F. (2013).
\newblock A continuous approximation model for time definite many-to-many transportation.
\newblock {\em Transportation Research Part B: Methodological}, 54:100--112.

\bibitem[Chen and Nie, 2017a]{chen2017analysis}
Chen, P.~W. and Nie, Y.~M. (2017a).
\newblock Analysis of an idealized system of demand adaptive paired-line hybrid transit.
\newblock {\em Transportation Research Part B: Methodological}, 102:38--54.

\bibitem[Chen and Nie, 2017b]{chen2017connecting}
Chen, P.~W. and Nie, Y.~M. (2017b).
\newblock Connecting e-hailing to mass transit platform: Analysis of relative spatial position.
\newblock {\em Transportation Research Part C: Emerging Technologies}, 77:444--461.

\bibitem[Chen et~al., 2012]{chen2012comparison}
Chen, Q., Lin, J., and Kawamura, K. (2012).
\newblock Comparison of urban cooperative delivery and direct delivery strategies.
\newblock {\em Transportation research record}, 2288(1):28--39.

\bibitem[Daganzo, 2005]{daganzo2005logistics}
Daganzo, C. (2005).
\newblock {\em Logistics systems analysis}.
\newblock Springer Science \& Business Media.

\bibitem[Daganzo and Newell, 1986a]{daganzo1986design}
Daganzo, C. and Newell, G. (1986a).
\newblock Design of multiple-vehicle delivery tours--i: a ring-radial network.
\newblock {\em Transportation Research Part B}, pages 345--364.

\bibitem[Daganzo, 1978]{daganzo1978approximate}
Daganzo, C.~F. (1978).
\newblock An approximate analytic model of many-to-many demand responsive transportation systems.
\newblock {\em Transportation Research}, 12(5):325--333.

\bibitem[Daganzo, 1984a]{daganzo1984distance}
Daganzo, C.~F. (1984a).
\newblock The distance traveled to visit n points with a maximum of c stops per vehicle: An analytic model and an application.
\newblock {\em Transportation science}, 18(4):331--350.

\bibitem[Daganzo, 1984b]{daganzo1984length}
Daganzo, C.~F. (1984b).
\newblock The length of tours in zones of different shapes.
\newblock {\em Transportation Research Part B: Methodological}, 18(2):135--145.

\bibitem[Daganzo, 1987a]{daganzo1987break}
Daganzo, C.~F. (1987a).
\newblock The break-bulk role of terminals in many-to-many logistic networks.
\newblock {\em Operations Research}, 35(4):543--555.

\bibitem[Daganzo, 1987b]{daganzo1987increasing}
Daganzo, C.~F. (1987b).
\newblock Increasing model precision can reduce accuracy.
\newblock {\em Transportation Science}, 21(2):100--105.

\bibitem[Daganzo et~al., 2012]{daganzo2012potential}
Daganzo, C.~F., Gayah, V.~V., and Gonzales, E.~J. (2012).
\newblock The potential of parsimonious models for understanding large scale transportation systems and answering big picture questions.
\newblock {\em EURO Journal on Transportation and Logistics}, 1(1-2):47--65.

\bibitem[Daganzo and Newell, 1986b]{daganzo1986configuration}
Daganzo, C.~F. and Newell, G.~F. (1986b).
\newblock Configuration of physical distribution networks.
\newblock {\em Networks}, 16(2):113--132.

\bibitem[del Castillo, 1999]{del1999heuristic}
del Castillo, J.~M. (1999).
\newblock A heuristic for the traveling salesman problem based on a continuous approximation.
\newblock {\em Transportation Research Part B: Methodological}, 33(2):123--152.

\bibitem[Flushberg~Sr and Wilson, 1976]{flushberg1976descriptive}
Flushberg~Sr, M. and Wilson, N.~H. (1976).
\newblock A descriptive supply model for demand-responsive transportation system planning.

\bibitem[Ghaffarinasab et~al., 2018]{ghaffarinasab2018continuous}
Ghaffarinasab, N., Van~Woensel, T., and Minner, S. (2018).
\newblock A continuous approximation approach to the planar hub location-routing problem: Modeling and solution algorithms.
\newblock {\em Computers \& Operations Research}, 100:140--154.

\bibitem[Gunes et~al., 2024]{gunes2024seattle}
Gunes, S., Fried, T., and Goodchild, A. (2024).
\newblock Seattle microhub delivery pilot: Evaluating emission impacts and stakeholder engagement.
\newblock {\em Case Studies on Transport Policy}, 15:101119.

\bibitem[Kaczmarski, 2024]{kaczmarski2024which}
Kaczmarski, M. (2024).
\newblock Which company is winning the restaurant food delivery war?
\newblock https://secondmeasure.com/datapoints/food-delivery-services-grubhub-uber-eats-doordash-postmates/.

\bibitem[Kawamura and Lu, 2007]{kawamura2007evaluation}
Kawamura, K. and Lu, Y. (2007).
\newblock Evaluation of delivery consolidation in us urban areas with logistics cost analysis.
\newblock {\em Transportation research record}, 2008(1):34--42.

\bibitem[Kohar and Jakhar, 2021]{kohar2021capacitated}
Kohar, A. and Jakhar, S.~K. (2021).
\newblock A capacitated multi pickup online food delivery problem with time windows: a branch-and-cut algorithm.
\newblock {\em Annals of Operations Research}, pages 1--22.

\bibitem[Larson and Odoni, 1981]{larson1981urban}
Larson, R.~C. and Odoni, A.~R. (1981).
\newblock {\em Urban operations research}.
\newblock Number Monograph.

\bibitem[Lin et~al., 2016]{lin2016sustainability}
Lin, J., Chen, Q., and Kawamura, K. (2016).
\newblock Sustainability si: logistics cost and environmental impact analyses of urban delivery consolidation strategies.
\newblock {\em Networks and Spatial Economics}, 16:227--253.

\bibitem[Liu et~al., 2021]{liu2021time}
Liu, S., He, L., and Max~Shen, Z.-J. (2021).
\newblock On-time last-mile delivery: Order assignment with travel-time predictors.
\newblock {\em Management Science}, 67(7):4095--4119.

\bibitem[Liu, 2019]{liu2019optimization}
Liu, Y. (2019).
\newblock An optimization-driven dynamic vehicle routing algorithm for on-demand meal delivery using drones.
\newblock {\em Computers \& Operations Research}, 111:1--20.

\bibitem[Naseraldin and Herer, 2011]{naseraldin2011location}
Naseraldin, H. and Herer, Y.~T. (2011).
\newblock A location-inventory model with lateral transshipments.
\newblock {\em Naval Research Logistics (NRL)}, 58(5):437--456.

\bibitem[Newell, 1971]{newell1971dispatching}
Newell, G.~F. (1971).
\newblock Dispatching policies for a transportation route.
\newblock {\em Transportation Science}, 5(1):91--105.

\bibitem[Newell and Daganzo, 1986]{newell1986design}
Newell, G.~F. and Daganzo, C.~F. (1986).
\newblock Design of multiple vehicle delivery tours—ii other metrics.
\newblock {\em Transportation Research Part B: Methodological}, 20(5):365--376.

\bibitem[Nourbakhsh and Ouyang, 2012]{nourbakhsh2012structured}
Nourbakhsh, S.~M. and Ouyang, Y. (2012).
\newblock A structured flexible transit system for low demand areas.
\newblock {\em Transportation Research Part B: Methodological}, 46(1):204--216.

\bibitem[Ouyang and Yang, 2023]{ouyang2023measurement}
Ouyang, Y. and Yang, H. (2023).
\newblock Measurement and mitigation of the ``wild goose chase'' phenomenon in taxi services.
\newblock {\em Transportation Research Part B: Methodological}, 167:217--234.

\bibitem[Pahwa and Jaller, 2022]{pahwa2022cost}
Pahwa, A. and Jaller, M. (2022).
\newblock A cost-based comparative analysis of different last-mile strategies for e-commerce delivery.
\newblock {\em Transportation Research Part E: Logistics and Transportation Review}, 164:102783.

\bibitem[Shi and Xu, 2024]{shi2024dine}
Shi, L. and Xu, Z. (2024).
\newblock Dine in or takeout? trends on restaurant service demand amid the covid-19 pandemic.
\newblock {\em Service Science}.

\bibitem[Shi et~al., 2022]{shi2022integer}
Shi, L., Xu, Z., Lejeune, M., and Luo, Q. (2022).
\newblock An integer l-shaped method for dynamic order fulfillment in autonomous last-mile delivery with demand uncertainty.
\newblock {\em arXiv preprint arXiv:2208.09067}.

\bibitem[Simoni and Winkenbach, 2023]{simoni2023crowdsourced}
Simoni, M.~D. and Winkenbach, M. (2023).
\newblock Crowdsourced on-demand food delivery: An order batching and assignment algorithm.
\newblock {\em Transportation Research Part C: Emerging Technologies}, 149:104055.

\bibitem[Smilowitz and Daganzo, 2007]{smilowitz2007continuum}
Smilowitz, K.~R. and Daganzo, C.~F. (2007).
\newblock Continuum approximation techniques for the design of integrated package distribution systems.
\newblock {\em Networks: An International Journal}, 50(3):183--196.

\bibitem[Steever et~al., 2019]{steever2019dynamic}
Steever, Z., Karwan, M., and Murray, C. (2019).
\newblock Dynamic courier routing for a food delivery service.
\newblock {\em Computers \& Operations Research}, 107:173--188.

\bibitem[Stein, 1978]{stein1978asymptotic}
Stein, D.~M. (1978).
\newblock An asymptotic, probabilistic analysis of a routing problem.
\newblock {\em Mathematics of Operations Research}, 3(2):89--101.

\bibitem[Toth and Vigo, 2002]{toth2002vehicle}
Toth, P. and Vigo, D. (2002).
\newblock {\em The vehicle routing problem}.
\newblock SIAM.

\bibitem[Tsao et~al., 2012]{tsao2012continuous}
Tsao, Y.-C., Mangotra, D., Lu, J.-C., and Dong, M. (2012).
\newblock A continuous approximation approach for the integrated facility-inventory allocation problem.
\newblock {\em European journal of operational research}, 222(2):216--228.

\bibitem[Wiles and van Brunt, 2001]{wiles2001optimal}
Wiles, P.~G. and van Brunt, B. (2001).
\newblock Optimal location of transshipment depots.
\newblock {\em Transportation Research Part A: Policy and Practice}, 35(8):745--771.

\bibitem[Wilson et~al., 1969]{wilson1969simulation}
Wilson, N., Sussman, J., Goodman, L., and Higonnet, B. (1969).
\newblock Simulation of a computer aided routing system (cars).
\newblock Technical report, Institute of Electrical and Electronics Engineers (IEEE).

\bibitem[Wilson et~al., 1976]{wilson1976advanced}
Wilson, N. H.~M., Weissberg, R.~W., and Hauser, J. (1976).
\newblock Advanced dial-a-ride algorithms research project.
\newblock Technical report.

\bibitem[Xie and Ouyang, 2015]{xie2015optimal}
Xie, W. and Ouyang, Y. (2015).
\newblock Optimal layout of transshipment facility locations on an infinite homogeneous plane.
\newblock {\em Transportation Research Part B: Methodological}, 75:74--88.

\bibitem[Xu et~al., 2020]{xu2020supply}
Xu, Z., Yin, Y., and Ye, J. (2020).
\newblock On the supply curve of ride-hailing systems.
\newblock {\em Transportation Research Part B: Methodological}, 132:29--43.

\bibitem[Ye et~al., 2024]{ye2024modeling}
Ye, A., Zhang, K., Chen, X.~M., Bell, M.~G., Lee, D.-H., and Hu, S. (2024).
\newblock Modeling and managing an on-demand meal delivery system with order bundling.
\newblock {\em Transportation Research Part E: Logistics and Transportation Review}, 187:103597.

\end{thebibliography}
\end{document}